\documentclass{imanum_NoJournalName}

\received{1 August 2014}

\usepackage{graphicx}
\usepackage{amsmath}
\usepackage{amsthm}
\usepackage{color}
\usepackage{float}
\usepackage{url}
\usepackage[dvips]{hyperref}

\newtheorem{thm}{Theorem}

\theoremstyle{definition}

\newtheorem{ex}{Example}
\newtheorem{rem}{Remark}

\numberwithin{equation}{section}

\def\E{\mathrm{E}}
\def\e{\mathrm{e}}
\def\i{\mathrm{i}}
\def\d{\mathrm{d}}

\allowdisplaybreaks[4]

\begin{document}

\title{A fast and accurate numerical method for the symmetric L\'evy processes 
based on the Fourier transform and sinc-Gauss sampling formula}
\shorttitle{A numerical method for the symmetric L\'evy processes}

\author{
{\sc Ken'ichiro Tanaka\thanks{Email: ketanaka@fun.ac.jp}} \\[2pt]
Future University Hakodate, 
116-2, Kamedanakano-cho, Hakodate, Hokkaido, Japan
}

\maketitle

\begin{abstract}
{In this paper, we propose a fast and accurate numerical method based on Fourier transform 
to solve Kolmogorov forward equations of symmetric scalar L\'evy processes. 
The method is based on the accurate numerical formulas for Fourier transform proposed by Ooura. 
These formulas are combined with nonuniform fast Fourier transform (FFT) and fractional FFT 
to speed up the numerical computations. 
Moreover, 
we propose a formula for numerical indefinite integration on equispaced grids
as a component of the method. 
The proposed integration formula is based on the sinc-Gauss sampling formula, 
which is a function approximation formula. 
This integration formula is also combined with the FFT. 
Therefore, all steps of the proposed method are executed using the FFT and its variants.
The proposed method allows us to be free from 
some special treatments for a non-smooth initial condition 
and numerical time integration. 
The numerical solutions obtained by the proposed method 
appeared to be exponentially convergent on the interval 
if the corresponding exact solutions do not have sharp cusps. 
Furthermore, the real computational times are approximately consistent with the theoretical estimates.  
}
{L\'evy process; Kolmogorov forward equation; nonuniform FFT; fractional FFT; sinc-Gauss sampling formula.}
\end{abstract}

\section{Introduction}

In this paper, we propose a fast and accurate numerical method based on the Fourier transform 
to solve the Kolmogorov forward equations of the symmetric scalar L\'evy processes. 
To propose the method, 
we use Ooura's accurate numerical formulas for the Fourier transform 
\citep{bib:OouraEuler2001, bib:OouraDE-FT2005}, 
and propose a numerical indefinite integration formula 
based on the sinc-Gauss sampling formula~\citep{bib:KTanaka_etal_SG_2008} 
to compute the integrals with respect to the L\'evy measures in the equation. 
Furthermore, we combine the fast Fourier transform (FFT) with these formulas 
to speed up the numerical computations.

A L\'evy process $\{ X_{t} \}_{t \geq 0}$ is essentially a stochastic process 
with stationary and independent increments~\citep{bib:Apple_LevyText_2009}, 
which is used to describe uncertain phenomena in various fields. 
Let us consider a brief non-rigorous review of 
the scalar L\'evy process $\{ X_{t} \}_{t \geq 0}$. 
The L\'evy-Khintchine theorem characterizes 
$ X_{t} $ by reals $b\in \mathbf{R}$, $a\geq 0$, 
and a Borel measure $\nu$ on $\mathbf{R} \setminus \{ 0 \}$ as
\begin{align}
\E[\e^{\i\, u\, X_{t}}] = \e^{t\, \psi(u)},
\end{align}
where $\psi$ is the characteristic exponent of $X_{1}$ defined by
\begin{align}
\psi(u) = \i\, b\, u - \frac{1}{2} a^{2} u^{2}
+ \int_{\mathbf{R} \setminus \{0\} } 
\left(
\e^{\i\, u\, y} - 1 - \i\, u\, y\, \mathbf{1}_{|y|\leq 1}(y) 
\right)\,
\nu(\d y).
\end{align}
Here, we assume that the measure $\nu$ satisfies
\begin{align}
\int_{\mathbf{R} \setminus \{0\} } 
\min\left\{ y^2,\, 1 \right\}
\, \nu(\d y) < \infty. 
\end{align}
Then, $\nu$ is called the L\'evy measure. 
An operator semigroup $\{ T_{t} \}_{t\geq 0}$ is associated with $\{ X_{t} \}_{t \geq 0}$, 
namely $(T_{t} f)(x) = \mathrm{E}[f(X_{t}) \mid X_{0} = x]$, 
where $f$ is a bounded continuous function on $\mathbf{R}$. 
The infinitesimal generator $A = \frac{\d}{\d t}(T_{t} f) \mid _{t=0}$ takes the form
\begin{align}
(Af)(x) = b\, f'(x) + a\, f''(x) + \int_{\mathbf{R} \setminus \{0\} } 
\left[
f(x+y) - f(x) - y\, \mathbf{1}_{|y|\leq 1}(y)\, f'(x) 
\right]\,
\nu(\d y).
\label{eq:LevyOP}
\end{align}
Then, the function $u(x, t) = (T_{t} f)(x)$ is 
the solution of the partial integro-differential equation (PIDE)
\begin{align}
& \frac{\partial u}{\partial t}(x, t) = A\, u(x, t)
\quad (x \in \mathbf{R},\ t \geq 0) 
\label{eq:Diffu_PIDE} 
\end{align}
with initial condition $u(x,0) = f(x)$.
Assuming that there exist the transition probability measure $p(x_{0},0; x, t)$ of $X_{t}$
with appropriate continuity and differentiability for each $t \geq 0$ 
and the adjoint operator $A^{\dagger}$ of $A$, 
we have
\begin{align}
\frac{\partial p}{\partial t}(x_{0},0; x, t) = A^{\dagger}\, p(x_{0},0; x, t)
\quad (x \in \mathbf{R},\ t \geq 0) 
\label{eq:Diffu_PIDE_Aad}
\end{align}
with initial condition $p(x_{0},0; x, 0) = \delta(x - x_{0})$, 
where $\delta$ is the Dirac delta function. 
Equation \eqref{eq:Diffu_PIDE_Aad} is the Kolmogorov forward equation, 
which is also known as the Fokker-Planck equation. 
Furthermore, the equation 
\begin{align}
\frac{\partial v}{\partial s} (x, t-s) = -A\, v(x,t-s)
\quad (x \in \mathbf{R},\ s \leq t)
\label{eq:Diffu_PIDE_A}
\end{align}
with initial condition $v(x,t) = f(x)$ is known as the Kolmogorov backward equation. 

In various fields such as physics, chemistry, biology, engineering, and economics, 
the Kolmogorov forward equations
including fractional derivatives are often considered 
to describe some unusual diffusion such as anomalous diffusion.  
For example, see 
\cite{bib:Gardiner_StocMeth_2009, 
bib:Kozu_etal_FLM_2006, 
bib:Lenzi_etal_FracFPE_2003, 
bib:Sabatier_etal_FracCal_2007, 
bib:Yan_FFPE_Laplace_2013} 
and the references therein.
Such equations are related to $\alpha$-stable processes belonging to the L\'evy processes. 
In finance, some L\'evy processes are used to describe the prices of risk assets.
Then, the Kolmogorov backward equations for these processes are considered 
as one of the useful methods for option pricing. 
See 
\cite{bib:ContVoltchkova_PIDELevy_2005, 
bib:GarreauKopriva2013, 
bib:Lee_etal_FastExpIntOP_2012}
and the references therein.  

In the fields mentioned above, 
many numerical methods for these Kolmogorov equations are studied. 
Popular examples of such methods are 
finite difference methods
\citep{bib:Gao_FPE_symLevy_2013, 
bib:Huang_FracLaplaceI_2013, 
bib:Li_FDM_FPE_2012, 
bib:Meerschaert_FDM_Flow_2004}, 
finite element methods
\citep{bib:Zhao_NumFPE_2012}, 
spectral methods 
\citep{bib:Bueno_FourierFrac_2014, 
bib:Huang_FDSpec_2014}, and 
other methods 
\citep{bib:Yan_FFPE_Laplace_2013}
for the forward equations describing anomalous diffusion etc. 
Further, such methods have been proposed 
for the backward equations in finance 
\citep{bib:LevyMattersI_2010, 
bib:GarreauKopriva2013, 
bib:Kwok_etal_2012, 
bib:Lee_etal_FastExpIntOP_2012}. 
In many of these methods, 
first, a time-evolution system of ordinary differential equations is derived from the given PIDE
by the discretization of the spatial variable with 
finite differences, finite elements, polynomial expansions, etc., 
and some quadrature formulas. 
Then, the system is numerically solved using some numerical time integration methods 
such as second-order finite difference methods. 
In addition, in the case 
where the closed form of the characteristic function of the L\'evy process can be obtained,  
methods based on the Fourier series or the Fourier transform are used in option pricing 
\citep{bib:CarrMadanPriceFFT1999, 
bib:ChourdakisPriceFracFFT2005, 
bib:FangOost_COS_2008, 
bib:Kwok_etal_2012}. 

In this paper, we propose a method based on the Fourier transform for equation~\eqref{eq:Diffu_PIDE_Aad}
for broader classes of the L\'evy measure $\nu$  
for which the closed form of the corresponding characteristic function may not be available. 
For simplicity, we consider the case $b = a = 0$ in~\eqref{eq:LevyOP} and measure $\nu$ has the form 
\begin{align}
\nu(\d y) = \frac{1}{|y|^{\gamma}}\, \mu(|y|)\, \d y, 
\label{eq:Lv_measure}
\end{align}
where $\gamma = 1$ or $\gamma = 2$ and $\mu \in L^{1}(0, \infty)$.
Then, the corresponding L\'evy process becomes symmetric. 
As examples of the L\'evy process with such measure, 
we can give  
the variance gamma (VG) process~\citep{bib:Apple_LevyText_2009}, 
also known as the symmetric Laplace motion~\citep{bib:Kozu_etal_FLM_2006}, 
with $\gamma = 1$ and $\mu(y) = \e^{- y}$, 
and the normal inverse Gaussian (NIG) process~\citep{bib:Apple_LevyText_2009} with $\gamma = 2$ and 
$\mu(y) = y\, K_{1}(y)/\pi$, where $K_{1}$ is the modified Bessel function of the second kind. 
Then, applying \eqref{eq:Lv_measure} to the operator $A$ in~\eqref{eq:LevyOP} 
and taking its adjoint, we have a special form of equation~\eqref{eq:Diffu_PIDE_Aad} as
\begin{align}
\frac{\partial p}{\partial t}(x, t) = E_{\gamma}^{+} p(x,t) + E_{\gamma}^{-} p(x,t)
\quad (x \in \mathbf{R},\ t \geq 0),
\label{eq:PIDE_LM}
\end{align}
where $p(x_{0},0; x,t)$ is denoted by $p(x,t)$ for conciseness, and
\begin{align}
E_{\gamma}^{\pm} q(x) = \int_{0}^{\infty} \frac{q(x \pm y) - q(x)}{y^{\gamma}}\, \mu(y)\, \d y. 
\label{eq:Op_E}
\end{align}
Note that the third term of the integrand in~\eqref{eq:LevyOP} vanishes 
because of the symmetry of $\nu$. 
Furthermore, we assume that $x_{0} = 0$ for simplicity
and $p(x,t) \to 0$ as $|x| \to 0$ for any $t \geq 0$ so that 
\begin{align}
\int_{-\infty}^{\infty} p(x,t)\, \d x = 1
\label{eq:TotalProb}
\end{align}
for any $t\geq 0$. 
Therefore, we consider equation \eqref{eq:PIDE_LM} with initial condition 
$p(x,0) = \delta(x)$ and auxiliary condition~\eqref{eq:TotalProb} in the rest of this paper.

There are two main objectives of the method proposed in this paper. 
The first is to show that 
a fast and accurate Fourier-based method can be realized for equation~\eqref{eq:PIDE_LM} 
including a (seemingly) singular integral in~\eqref{eq:Op_E}. 
Using the Fourier transform, 
we do not need specific requirements for 
the non-smooth initial condition $p(x, 0) = \delta(x)$, 
whereas some of the existing methods with time integration cited above 
need a priori artificial approximation of the solution $p(x, t)$ for a small time $t$. 
In addition, since equation~\eqref{eq:PIDE_LM} is linear and contains only constant coefficients, 
which is also the case for the general $A$ in~\eqref{eq:LevyOP}, 
the Fourier-based method does not need numerical time integration. 
Therefore, we can compute approximate solutions of $p(x, t)$ 
for any $t$ with the same computational cost, 
and errors of numerical time integration do not occur. 
The second objective of the proposed method is to 
present new applications of Ooura's methods~\citep{bib:OouraEuler2001, bib:OouraDE-FT2005} 
for the Fourier transforms to obtain the solution of PIDEs. 
The high accuracy of the proposed method is due to the very fast convergence of Ooura's methods, 
and speeding up of computations by Ooura's methods is realized by 
combining them with 
the nonuniform FFT~\citep{bib:DuttRokhlin_NFFT_1993, 
bib:DuttRokhlin_NFFT_1995, 
bib:GreengardLee_NFFT_2004, 
bib:Potts_etal_NFFT_2001, 
bib:Steidl_NFFT_1998}
or
the fractional FFT~\citep{bib:BailSwarFRFT1991, 
bib:ChourdakisPriceFracFFT2005, 
bib:KTanaka_EulerFFT_2014}. 
In addition, in order to treat a (seemingly) singular integral in~\eqref{eq:Op_E}, 
we propose an indefinite integration formula 
using the sinc-Gauss sampling formula~\citep{bib:KTanaka_etal_SG_2008}, 
which is also accurate and combined with the FFT. 
Thus, as shown precisely in Sections~\ref{sec:Outline} and~\ref{sec:NumScheme}, 
all steps of the proposed method are executed by the FFT and its variants.

The remainder of this paper is organized as follows.
Section~\ref{sec:Outline} contains the outline of the proposed method, 
which consists of three steps. 
Section~\ref{sec:NumScheme} details the three steps. 
Section~\ref{sec:NumEx} shows 
the actual performance of the proposed method through numerical examples, 
and Section~\ref{sec:concl} concludes this paper.

\section{Outline of the proposed method}
\label{sec:Outline}

In order to obtain a fast and accurate numerical method to solve~\eqref{eq:PIDE_LM}, 
we considered a method based on the Fourier transform for using the FFT. 
First, we derive the formula for the solution of \eqref{eq:PIDE_LM} using the Fourier transform 
\begin{align}
[\mathcal{F} f](\omega) = \int_{-\infty}^{\infty} f(x)\, \e^{-\i\, \omega\, x}\, \d x. 
\end{align}
Taking the Fourier transform for both sides of \eqref{eq:PIDE_LM} 
with respect to the spatial variable $x$, 
we have 
\begin{align}
\frac{\partial\, [\mathcal{F} p]}{\partial t}(\omega, t)
=
[\mathcal{G}_{\gamma} \mu](\omega)\,  
[\mathcal{F} p](\omega, t), 
\label{eq:PIDE_LM_FT_Gen}
\end{align}
where 
\begin{align}
[\mathcal{G}_{\gamma} \mu](\omega) 
= 
\int_{0}^{\infty}
\frac{\e^{-\i\, \omega\, y} - 1}{y^{\gamma}}\, \mu(y)\, \d y
+
\int_{0}^{\infty}
\frac{\e^{+\i\, \omega\, y} - 1}{y^{\gamma}}\, \mu(y)\, \d y.
\label{eq:PsiMu}
\end{align}
For $\gamma = 1$, noting that $\mu \in L^{1}(0, \infty)$, we have 
\begin{align}
\int_{0}^{\infty}
\frac{\e^{\mp \i\, \omega\, y} - 1}{y}\, \mu(y)\, \d y
&=
\int_{0}^{\infty}
\left(
\mp \i 
\int_{0}^{\omega} \e^{\mp \i\, \zeta \, y}\, \d \zeta 
\right)
\, \mu(y)\, \d y \notag \\
& = 
\mp \i 
\int_{0}^{\omega} 
\left(
\int_{0}^{\infty} \mu(y)\, \e^{\mp \i\, \zeta \, y}\, \d y 
\right)
\d \zeta, 
\label{eq:TransInt}
\end{align}
and therefore,
\begin{align}
[\mathcal{G}_{1} \mu](\omega) 
=
2\, \mathop{\mathrm{Im}}
\int_{0}^{\omega} 
\left(
\int_{0}^{\infty} \mu(y)\, \e^{- \i\, \zeta \, y}\, \d y 
\right)
\d \zeta.
\label{eq:G_alpha_1}
\end{align}
For $\gamma = 2$, noting that 
\begin{align}
& \int_{0}^{\infty}
\frac{\e^{- \i\, \omega\, y} - 1}{y^2}\, \mu(y)\, \d y
+
\int_{0}^{\infty}
\frac{\e^{+ \i\, \omega\, y} - 1}{y^2}\, \mu(y)\, \d y \notag \\
& =
\int_{0}^{\infty}
\frac{\e^{- \i\, \omega\, y} - 1 + \i\, \omega\, y}{y^2}\, \mu(y)\, \d y
+
\int_{0}^{\infty}
\frac{\e^{+ \i\, \omega\, y} - 1 - \i\, \omega\, y}{y^2}\, \mu(y)\, \d y \notag 
\end{align}
and $\mu \in L^{1}(0, \infty)$, we have 
\begin{align}
\int_{0}^{\infty}
\frac{\e^{\mp \i\, \omega\, y} - 1 \pm \i\, \omega\, y}{y^2}\, \mu(y)\, \d y
&=
\int_{0}^{\infty}
\left(
- \int_{0}^{\omega} \int_{0}^{\eta} \e^{\mp \i\, \zeta \, y}\, \d \zeta\, \d \eta 
\right)
\mu(y)\, \d y \notag \\
& = 
-
\int_{0}^{\omega} \int_{0}^{\eta} 
\left(
\int_{0}^{\infty} \mu(y)\, \e^{\mp \i\, \zeta \, y}\, \d y 
\right)
\d \zeta\, \d \eta. 
\label{eq:TransInt_2}
\end{align}
Therefore, we have 
\begin{align}
[\mathcal{G}_{2} \mu](\omega) 
=
-2\, \mathop{\mathrm{Re}}
\int_{0}^{\omega} \int_{0}^{\eta} 
\left(
\int_{0}^{\infty} \mu(y)\, \e^{- \i\, \zeta \, y}\, \d y 
\right)
\d \zeta\, \d \eta. 
\label{eq:G_alpha_2}
\end{align}
Using expression~\eqref{eq:G_alpha_1} or~\eqref{eq:G_alpha_2}, 
we can derive the form of the solution $p(x,t)$ 
from~\eqref{eq:PIDE_LM_FT_Gen} as 
\begin{align}
p(x,t) 
=
\mathcal{F}^{-1}
\left[\, \exp 
\left( t\, [\mathcal{G}_{\gamma} \mu]
\right)\, 
\right] (x)
=
\frac{1}{2\pi} \int_{-\infty}^{\infty} 
\exp \left( t\, [\mathcal{G}_{\gamma} \mu](\omega) \right)\, \e^{\i\, x\, \omega}\, \d \omega.
\label{eq:ExprExactSol}
\end{align}

Then, we can consider the following numerical method 
to obtain an approximation of the solution~\eqref{eq:ExprExactSol}.
Let $N>0$, $\hat{h} > 0$, and $t > 0$ be 
an integer, a grid spacing, and a time, respectively. 
Suppose that we need to compute approximate values of $p(x,t)$ for 
$x = n \hat{h}\ (n = -N+1,\ldots , N)$. 
\begin{description}
\setlength{\parskip}{0pt}
\setlength{\itemsep}{3pt}
\item[Step 1] Computation of the Fourier transform
\begin{align}
\int_{0}^{\infty} \mu(y)\, \e^{- \i\, \zeta \, y}\, \d y
\label{eq:first_FT}
\end{align}
in~\eqref{eq:G_alpha_1} or~\eqref{eq:G_alpha_2}. 
Choose a grid spacing $\tilde{h} > 0$ and integers $M_{-}, M_{+}, N_{\gamma} > 0$. 
Then, use {\it the double exponential (DE) formula for the Fourier transforms}~\citep{bib:OouraDE-FT2005} 
with sampling points 
$y = y_{j}\ (j = -M_{-}, \ldots, M_{+}-1)$ 
and {\it the nonuniform FFT}~\citep{bib:DuttRokhlin_NFFT_1993, 
bib:DuttRokhlin_NFFT_1995, 
bib:GreengardLee_NFFT_2004, 
bib:Potts_etal_NFFT_2001, 
bib:Steidl_NFFT_1998}
to obtain approximate values of~\eqref{eq:first_FT} 
for $\zeta = k \tilde{h}\ (k = -N_{\gamma}+1,\ldots, N_{\gamma})$. 
The definitions of $\tilde{h}, M_{-}, M_{+}, N_{\gamma}$, and $y_{j}$ are presented in 
Sections~\ref{sec:Step1} and~\ref{sec:SGII}.
\item[Step 2] 
Computation of the indefinite integral of \eqref{eq:first_FT}
in~\eqref{eq:G_alpha_1} or~\eqref{eq:G_alpha_2}. 
Use the computed values in Step~1 and 
{\it an indefinite integration by the sinc-Gauss sampling formula}
proposed in Section~\ref{sec:SGII} to obtain 
approximate values of $[\mathcal{G}_{\gamma} \mu](\omega)$  for 
$\omega = \ell \tilde{h}\ (\ell = -N_{\gamma}/2^{\gamma}+1,\ldots, N_{\gamma}/2^{\gamma})$, 
where we suppose that $N_{\gamma}$ can be divided by $2^{\gamma}$. 
Note that the approximate values on the equispaced grid are 
obtained from the approximate values of the integrand on the same grid. 
The definition of $N_{\gamma}$ is presented in Section~\ref{sec:SGII}. 
\item[Step 3]
Computation of the inverse Fourier transform~\eqref{eq:ExprExactSol}.
Use the computed values in Step~2, 
{\it the formula for the Fourier transform with continuous Euler transform}~\citep{bib:OouraEuler2001}, 
and {\it the fractional FFT}~\citep{bib:BailSwarFRFT1991, 
bib:ChourdakisPriceFracFFT2005, 
bib:KTanaka_EulerFFT_2014}
to obtain the approximate values of the solution $p(x,t)$
for $x = n \hat{h}\ (n = -N+1,\ldots , N)$. 
\end{description}
In summary, the proposed method is illustrated by the diagram below. 
The details of the three steps are shown in Section~\ref{sec:NumScheme}. 

\newpage

\begin{center}

\begin{minipage}[c]{.35\linewidth}
The Kolmogorov forward equation~\eqref{eq:PIDE_LM} 
for the L\'evy process with 
measure given by~\eqref{eq:Lv_measure}
\end{minipage}
\begin{minipage}[c]{.25\linewidth}
\begin{center}
$ \xrightarrow{
\begin{subarray}{l}
\text{\bf (Step 1)}\\
\text{\bf DE formula for FT} \\
\text{\bf $+$ nonuniform FFT}
\end{subarray}
} $
\end{center}
\end{minipage}
\begin{minipage}[c]{.35\linewidth}
Approximate values of \eqref{eq:first_FT}: 
\[
\int_{0}^{\infty} \mu(y)\, \e^{- \i\, \zeta \, y}\, \d y
\]
\end{minipage} 

\vspace{10pt}

\begin{minipage}[c]{.35\linewidth}
\ 
\end{minipage}
\begin{minipage}[c]{.25\linewidth}
\ 
\end{minipage}
\begin{minipage}[c]{.35\linewidth}
\rotatebox[origin = c]{90}{$\longleftarrow$}\  
\begin{minipage}[c]{\linewidth}
{\scriptsize \bf (Step 2) \\
sinc-Gauss indefinite integration}
\end{minipage}
\end{minipage} \\

\vspace{10pt}

\begin{minipage}[c]{.35\linewidth}
Approximate values of the solution~\eqref{eq:ExprExactSol}: 
\[
\displaystyle p(x, t) = \mathcal{F}^{-1}
\left[\, \exp 
\left( t\, [\mathcal{G}_{\gamma} \mu]
\right)\, 
\right] (x)
\]
\end{minipage}
\begin{minipage}[c]{.25\linewidth}
\begin{center}
$ \xleftarrow{
\begin{subarray}{l}
\text{\bf (Step 3)} \\
\text{\bf Formula with continuous} \\
\text{\bf Euler transform for FT} \\ 
\text{\bf $+$ fractional FFT}
\end{subarray}
} $
\end{center}
\end{minipage}
\begin{minipage}[c]{.35\linewidth}
Approximate values of 
$\displaystyle 
[\mathcal{G}_{\gamma} \mu](\omega) $
in~\eqref{eq:G_alpha_1} or~\eqref{eq:G_alpha_2}. 
\end{minipage}

\end{center}

\section{Proposed method}
\label{sec:NumScheme}

\subsection{Step 1: Computation of the Fourier transform~\eqref{eq:first_FT} in~\eqref{eq:G_alpha_1} or~\eqref{eq:G_alpha_2}}
\label{sec:Step1}

The DE formula for the Fourier transforms and
the nonuniform FFT for the Fourier transform~\eqref{eq:first_FT}
are described in Sections~\ref{sec:DE-FT} and~\ref{sec:NFFT}, respectively. 
The contribution of this paper is speeding up the computation 
through the use of the DE formula 
by combining it with nonuniform FFT.  

\subsubsection{DE formula for the Fourier transforms by Ooura}
\label{sec:DE-FT}

We begin with the review of the DE formula for the Fourier transforms~\eqref{eq:first_FT}
proposed by \cite{bib:OouraDE-FT2005}.
Let $\zeta_{0}$ and $h$ be positive constants and let 
the function $\varphi: \mathbf{R} \to (0,\infty)$ be defined by
\begin{align}
\varphi(t) 
=
\frac{t}{1 - \exp ( -2t - \alpha (1-\e^{-t}) - \beta (\e^{t} - 1) )}, 
\label{eq:DEFT_trans}
\end{align}
where $\beta = 0.25$ and 
\begin{align}
\alpha = \frac{\beta}{\sqrt{1 + \log(1 + \pi/(\zeta_{0} h))/(4 \zeta_{0} h)}}.
\label{eq:def_DEFT_alpha}
\end{align}
Then, the following formula approximates the integral~\eqref{eq:first_FT} 
for $\zeta \in (\delta_{1}, 2 \zeta_{0} - \delta_{2})$ for some $\delta_{1}, \delta_{2} \geq 0$:
\begin{align}
& \int_{0}^{\infty} \mu(y)\, \e^{- \i\, \zeta \, y}\, \d y \notag \\
& \approx
-\frac{2\pi \i}{\zeta_{0}}
\sum_{j = -M_{-}}^{M_{+}-1} 
\left[
\mu\left( \frac{\pi}{\zeta_{0} h} \varphi(j h) \right)\, 
\sin\left( \frac{\pi}{2 h} \hat{\varphi}(j h) \right)\, 
\varphi'(j h)\, 
\exp\left( \frac{\pi\, \i}{2 h} \hat{\varphi}(j h) \right)
\right]
\exp\left( -\frac{\pi\, \i\, \zeta}{\zeta_{0} h} \varphi(j h) \right).
\label{eq:DE_FT}
\end{align}
where $\hat{\varphi}(t) = \varphi(t) - t$.
The integers $M_{-}$ and $M_{+}$ are determined in an appropriate manner. 
Formula~\eqref{eq:DE_FT} is 
the DE formula for the Fourier transforms, which is derived as follows. 
First, applying the variable transformation 
$y = (\pi/(\zeta_{0}h))\, \varphi(t)$ to integral~\eqref{eq:first_FT}, 
we have 
\begin{align}
\int_{0}^{\infty} \mu(y)\, \e^{-\i\, \zeta \, y}\, \d y
=
\int_{-\infty}^{\infty} 
\mu\left( P \varphi(t) \right)\, 
\exp\left( -\i\, \zeta\, P \varphi(t) \right)
P \varphi'(t)\, \d t, 
\label{eq:DE_trans_1}
\end{align}
where $P = \pi/(\zeta_{0}h)$. 
Let $\mathcal{M}(\zeta)$ denote expression~\eqref{eq:DE_trans_1}. 
Next, from $\mathcal{M}(\zeta)$ subtract 
\begin{align}
\mathcal{N}(\zeta) = 
\int_{-\infty}^{\infty} 
\mu\left( P \varphi(t) \right)\, 
\exp\left( -\i\, \zeta\, P \varphi(t) + \i\, \zeta_{0}\, P \hat{\varphi}(t) \right)
P \varphi'(t)\, \d t, 
\end{align}
which is very small for $\zeta \in (\delta_{1}, 2\zeta_{0} - \delta_{2})$ and a large $P$. 
Then, discretizing 
\begin{align}
\mathcal{M}(\zeta) - \mathcal{N}(\zeta)
=
-\frac{2\pi \i}{\zeta_{0} h}
\int_{-\infty}^{\infty} 
\mu\left( \frac{\pi}{\zeta_{0} h} \varphi(t) \right)\, 
\sin\left( \frac{\pi}{2 h} \hat{\varphi}(t) \right)\, 
\varphi'(t)\, 
\exp\left( -\frac{\pi\, \i\, \zeta}{\zeta_{0} h} \varphi(t) + \frac{\pi\, \i}{2 h} \hat{\varphi}(t) \right)
\, \d t
\label{eq:before_disc}
\end{align}
by the mid-point rule with grid spacing $h$, 
we have~\eqref{eq:DE_FT}. 
Since $\hat{\varphi}(t) \to 0$ as $t \to -\infty$ and 
$\varphi'(t) \to 0$ as $t \to +\infty$, 
the factor
\(
\left[
\sin\left( \frac{\pi}{2 h} \hat{\varphi}(t) \right) 
\varphi'(t)
\right]
\)
in~\eqref{eq:before_disc} converges rapidly 
(``double exponentially'')
to $0$ as $t \to \pm \infty$. 
Therefore, the discretization of~\eqref{eq:before_disc} by the mid-point rule 
can yield accurate approximation~\eqref{eq:DE_FT}
for some $h$ independent of $\zeta$, and sufficiently large $M_{+}$ and $M_{-}$.  
In \cite{bib:OouraDE-FT2005}, 
the error of approximation~\eqref{eq:DE_FT} is bounded by 
$c_{0}' \e^{-c_{0}/h} + c_{1}' \e^{-c_{1} \zeta/h} + c_{2}' \e^{-c_{2}(2\zeta_{0} - \zeta)/h}$
for some $c_{i}, c_{i}'$ depending on $\mu$, 
and it is illustrated by some numerical examples. 
A theoretically rigorous analysis for the error, 
however, is not described in~\cite{bib:OouraDE-FT2005}. 

\begin{figure}[ht]
\begin{center}

\begin{minipage}[t]{0.45\linewidth}
\begin{center}
\includegraphics[width = \linewidth]{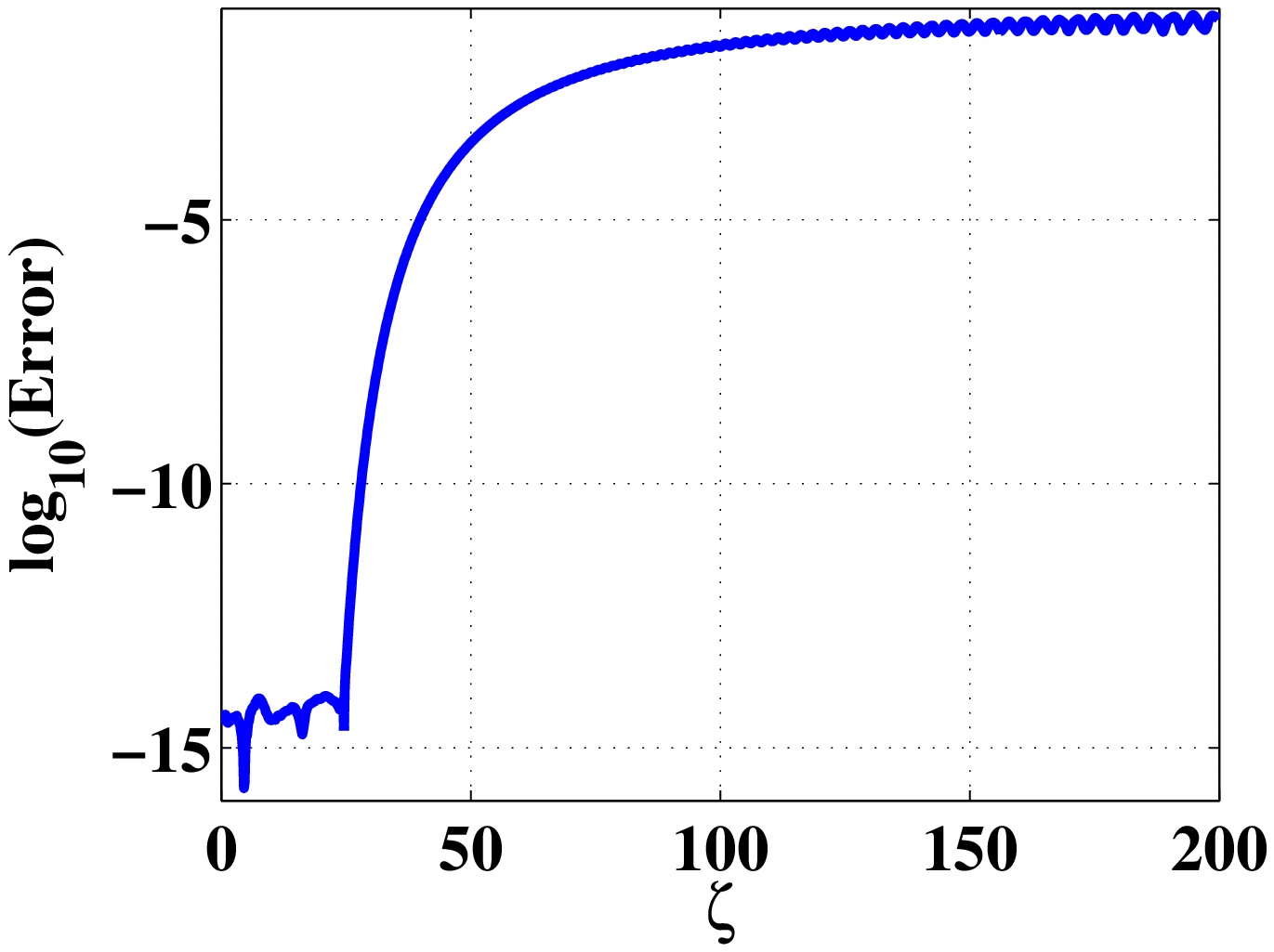}\\
(a) $\zeta_{0} = N_{\gamma} \tilde{h} / 15 \approx 10.0$
\end{center}
\end{minipage}\quad 
\begin{minipage}[t]{0.45\linewidth}
\begin{center}
\includegraphics[width = \linewidth]{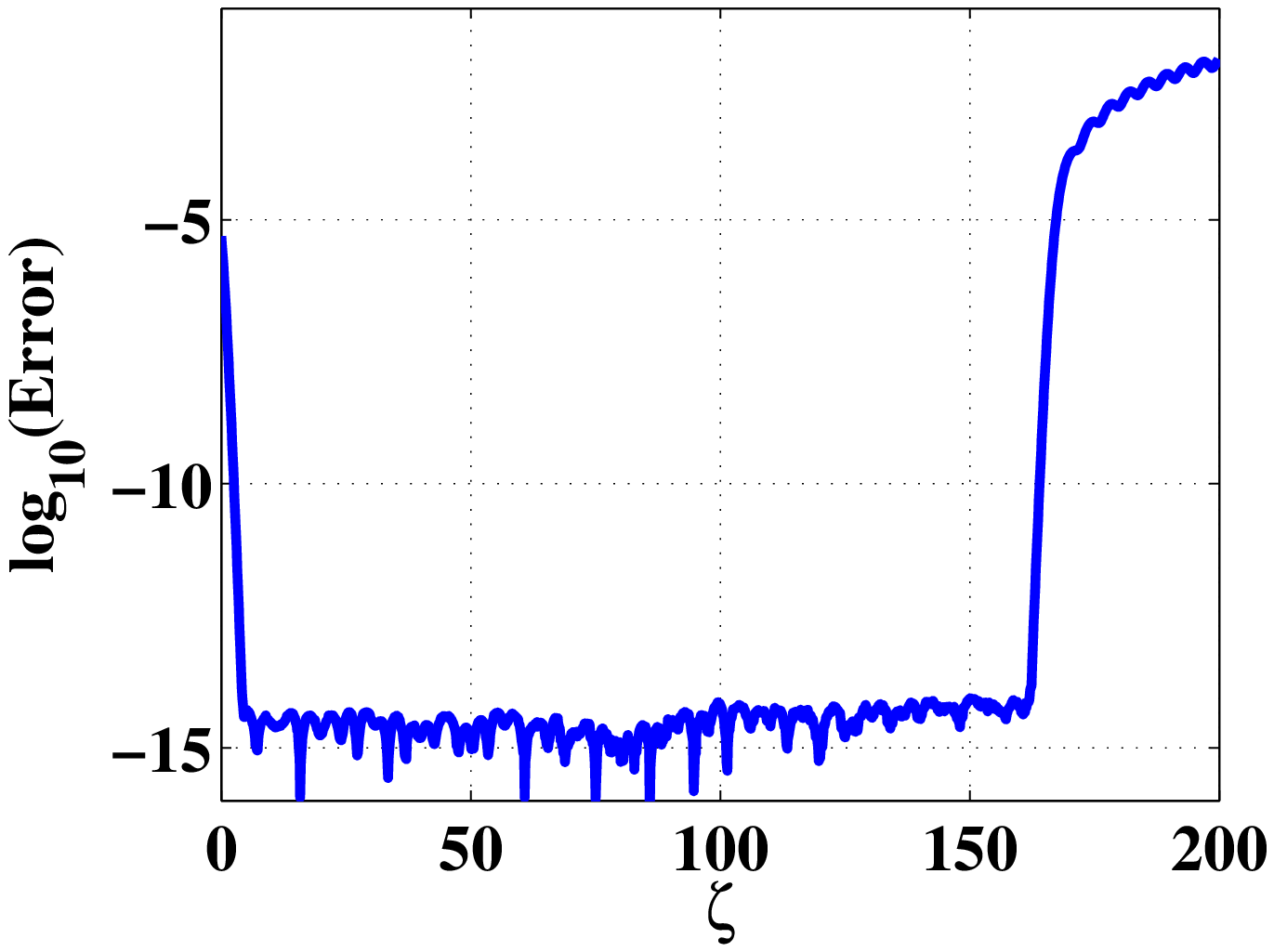}\\
(b) $\zeta_{0} = N_{\gamma} \tilde{h} / 1.8 \approx 83.4$
\end{center}
\end{minipage}
\vspace{5pt}
\caption{
Errors of the approximate values~\eqref{eq:DE_FT} 
for the Fourier transform~\eqref{eq:first_FT} of $\mu(y) = \e^{-y}$. 
The parameters $M_{-}, M_{+}$, and $h$ in~\eqref{eq:DE_FT} are defined by 
$M_{-} = M_{+} = 2^{10}$, and $h = \log(10^{3} M)/M \approx 0.007$, where $M = M_{-} + M_{+} = 2^{11}$. 
In addition, two cases for $\zeta_{0}$ in~\eqref{eq:DE_FT} are considered: 
(a) $\zeta_{0} = N_{\gamma} \tilde{h} / 15 \approx 10.0$, and 
(b) $\zeta_{0} = N_{\gamma} \tilde{h} / 1.8 \approx 83.4$, 
where $N_{\gamma} = M/2 = 2^{9}$ and $\tilde{h} = \sqrt{14 \pi / M} \approx 0.147$. 
The settings of $\zeta_{0}$ in (a) and (b) correspond 
to~\eqref{eq:set_zeta0_1} and~\eqref{eq:set_zeta0_2}, respectively. 
The function $\mu$ and the parameters above are also adopted 
in Example~\ref{ex:VG} in Section~\ref{sec:NumEx}. 
The approximate values~\eqref{eq:DE_FT} and their errors are computed for 
$\zeta = 0, \tilde{h}, \ldots, N_{\gamma} \tilde{h}, \ldots, M \tilde{h}$ with double precision
and shown for $\zeta = 0, \tilde{h}, \ldots, N_{\gamma} \tilde{h},\ldots , \lfloor M/1.5 \rfloor \tilde{h}$
in the graphs above. 
In case (a), the errors for $\zeta \in [0, 2\zeta_{0}]$ are small  
whereas they rapidly become large as $\zeta$ increases. 
In case (b), almost all errors for $\zeta \in [0, 2\zeta_{0}]$ are small 
except for some $\zeta$ near both sides of $[0, 2\zeta_{0}]$. 
Note that $N_{\gamma} \tilde{h} \approx 150$.
}
\label{fig:err_DEFT}
\end{center}
\end{figure}

Then, noting that
\begin{align}
\int_{0}^{\infty} \mu(y)\, \e^{- \i\, (-\zeta) \, y}\, \d y 
=
\overline{
\int_{0}^{\infty} \mu(y)\, \e^{- \i\, \zeta \, y}\, \d y }, 
\label{eq:FT_conjugate}
\end{align}
we can achieve Step~1 by computing the values of~\eqref{eq:DE_FT} for 
$\zeta = k \tilde{h}\ (k=0,\ldots, N_{\gamma})$ 
and taking their complex conjugates
for $\zeta = -k \tilde{h}\ (k=0,\ldots, N_{\gamma})$. 
In computing the values of~\eqref{eq:DE_FT}, 
we need to choose $\zeta_{0}$ 
so that $k \tilde{h} \in (\delta_{1}, 2\zeta_{0} - \delta_{2}) $ for $k=0,\ldots, N_{\gamma}$. 
The possible values of the nonnegative reals $\delta_{1}$ and $\delta_{2}$, 
however, are not theoretically estimated. 
According to some numerical examples including those in~\cite{bib:OouraDE-FT2005}, 
when $\zeta_{0}$ is small, $\delta_{1}$ can be taken as $\delta_{1} = 0$ and $\delta_{2}$ can be small. 
As $\zeta_{0}$ becomes large, unfortunately, $\delta_{1}$ and $\delta_{2}$ need to be large. 
These facts are illustrated by Figure~\ref{fig:err_DEFT}. 
Therefore, 
if we let $\zeta_{0}$ be a single value near to $N_{\gamma} \tilde{h} / 2$ when $N_{\gamma} \tilde{h}$ is large, 
we cannot have accurate approximations of~\eqref{eq:first_FT} for 
$\zeta = k \tilde{h}$ for $k$'s near to $0$ or $N_{\gamma}$. 
Then, we use 
\begin{align}
& \zeta_{0} = N_{\gamma} \tilde{h} / 15 
\text{\quad to compute~\eqref{eq:first_FT} for \quad} 
\zeta = 0, \tilde{h}, \ldots, \lfloor N_{\gamma} /8 \rfloor \tilde{h}, \label{eq:set_zeta0_1} \\
& \zeta_{0} = N_{\gamma} \tilde{h} / 1.8 
\text{\quad to compute~\eqref{eq:first_FT} for \quad} 
\zeta = (\lfloor N_{\gamma} /8 \rfloor  + 1) \tilde{h}, \ldots, N_{\gamma} \tilde{h}. \label{eq:set_zeta0_2} 
\end{align}
Figure~\ref{fig:err_DEFT} also illustrates these settings, 
which are experientially determined and not based on theoretical criteria.

Note that the naive computation of~\eqref{eq:DE_FT} 
for~\eqref{eq:set_zeta0_1} and~\eqref{eq:set_zeta0_2} 
requires $\mathrm{O}(N_{\gamma}^{2})$ operations if $M_{+} + M_{-} \propto N_{\gamma}$. 
Then, what remains in Step 1 is to speed up the numerical computation. 
Thus, we use the technique of the nonuniform FFT explained 
in Section~\ref{sec:NFFT} below.

\subsubsection{Nonuniform FFT}
\label{sec:NFFT}

\begin{figure}[ht]
\begin{center}

\begin{minipage}[t]{0.45\linewidth}
\begin{center}
\includegraphics[width = \linewidth]{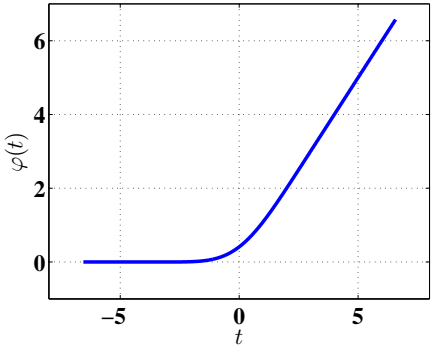}\\
(a) Example of the function $\varphi$ in~\eqref{eq:DEFT_trans}.
\end{center}
\end{minipage}\quad
\begin{minipage}[t]{0.45\linewidth}
\begin{center}
\includegraphics[width = \linewidth]{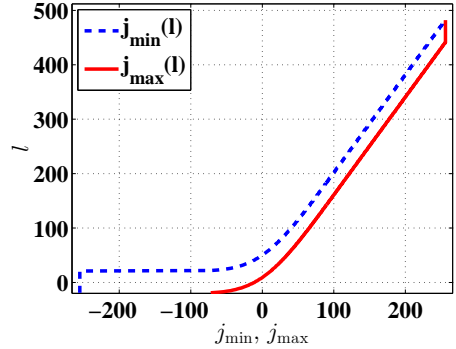}\\
(b) Example of the index set $\tilde{J}(l)$ in~\eqref{eq:aug_J_l}. 
\end{center}
\end{minipage}
\end{center}

\caption{
Examples of $\varphi$ in~\eqref{eq:DEFT_trans} and $\tilde{J}(l)$ in~\eqref{eq:aug_J_l}. 
(a) 
To fix $\varphi$, the following parameters are used:
$M_{-} = M_{+} = 2^{8}$, 
$M = M_{-} + M_{+}$, 
$N_{\gamma} = M/2$, 
$h = \log(10^{3} M)/M \approx 0.026$, 
$\zeta_{0} = \text{(the value of~\eqref{eq:set_zeta0_2})} \approx 41.7$, 
$\beta = 0.25$, and 
$\alpha = \text{(the value of~\eqref{eq:def_DEFT_alpha})} \approx 0.19$.
(b) 
Functions $j_{\min}$ and $j_{\max}$ are defined 
by~\eqref{eq:j_min} and~\eqref{eq:j_max}, respectively. 
To fix them, 
$\tilde{h}$ is determined by 
$\tilde{h} = \sqrt{14 \pi / M} \approx 0.29$, 
and $a$ and $\check{h}$ are determined by~\eqref{eq:NFFT_para_1}:  
$a = 2\pi/M \approx 0.012$, $\check{h} = 1$.
Finally, $b=20$.
The parameters stated above are also adopted in Example~\ref{ex:VG} in Section~\ref{sec:NumEx}. 
}
\label{fig:DE_and_IndexSetJ}
\end{figure}

Let the sum in \eqref{eq:DE_FT} for $\zeta = k \tilde{h}$ be rewritten as 
\begin{align}
\hat{\mu}_{k}
=
\sum_{j = -M_{-}}^{M_{+}-1} 
\Phi[\mu]_{j}\, 
\exp\left( -\i\, k \tilde{h}\, y_{j} \right), 
\label{eq:DE_FT_rewrite}
\end{align}
where 
\begin{align}
\Phi[\mu]_{j} 
& = 
-\frac{2\pi \i}{\zeta_{0}}
\left[
\mu\left( \frac{\pi}{\zeta_{0} h} \varphi(j h) \right)\, 
\sin\left( \frac{\pi}{2 h} \hat{\varphi}(j h) \right)\, 
\varphi'(j h)\, 
\exp\left( \frac{\pi\, \i}{2 h} \hat{\varphi}(j h) \right)
\right], \\
y_{j} 
& =
\frac{\pi}{\zeta_{0} h} \varphi(j h).  \label{eq:y_j}
\end{align}
The {\it nonuniform FFT}~\citep{bib:DuttRokhlin_NFFT_1993, 
bib:DuttRokhlin_NFFT_1995, 
bib:GreengardLee_NFFT_2004, 
bib:Potts_etal_NFFT_2001,
bib:Steidl_NFFT_1998}
is a fast method 
to compute the DFT  in \eqref{eq:DE_FT_rewrite} 
with a nonuniform grid such as $\{ y_{j} \}$ defined by \eqref{eq:y_j}. 
To use the technique of the nonuniform FFT, 
setting
\begin{align}
\tilde{\Phi}[\mu]_{j} 
& = 
\Phi[\mu]_{j}\, \exp \left( -\i \, \lfloor N_{\gamma}/2 \rfloor\, \tilde{h}\, y_{j} \right), \\
k' 
&=
k - \lfloor N_{\gamma}/2 \rfloor, \label{eq:k_prime}
\end{align}
we note the following relation: 
\begin{align}
\sqrt{\frac{\tau}{\pi}}\, \e^{-\tau\, (a k')^2}\, 
\hat{\mu}_{k}
& = 
\sum_{j = -M_{-}}^{M_{+}-1} 
\tilde{\Phi}[\mu]_{j}\, 
\sqrt{\frac{\tau}{\pi}}\, 
\exp\left( -\tau\, (a k')^2 - \i\, (a k')\, (\tilde{h} y_{j}/a) \,   \right)
\notag \\
& =
\sum_{j = -M_{-}}^{M_{+}-1} 
\tilde{\Phi}[\mu]_{j}\, 
\frac{1}{2\pi} \int_{-\infty}^{\infty} 
\exp \left( -\frac{(v - \tilde{h} y_{j}/a)^2}{4\tau} \right)\, 
\exp \left( - \i \, a k'\, v \right)\, \d v 
\notag \\
& =
\frac{1}{2\pi} \int_{-\infty}^{\infty} 
\left[
\sum_{j = -M_{-}}^{M_{+}-1} 
\tilde{\Phi}[\mu]_{j}\, 
\exp \left( -\frac{(v - \tilde{h} y_{j}/a)^2}{4\tau} \right)
\right]
\exp \left( - \i \, a k'\, v \right)\, \d v
\notag \\
& \approx
\frac{\check{h}}{2\pi} 
\sum_{l = -L_{-}}^{L_{+}}
\left[
\sum_{j = -M_{-}}^{M_{+}-1} 
\tilde{\Phi}[\mu]_{j}\, 
\exp \left( -\frac{(l \check{h} - \tilde{h} y_{j}/a)^2}{4\tau} \right)
\right]
\exp \left( - \i \, a \check{h}\, k' l \right), 
\label{eq:NfftKeyRel}
\end{align}
where $\tau$, $a$, $\check{h}$, and $L_{\pm}$
are positive constants to be determined appropriately. 
To decrease the computational cost of 
the sum in $[\ ]$ in \eqref{eq:NfftKeyRel}, 
we neglect the sufficiently small summands present in it
\begin{align}
\sum_{j = -M_{-}}^{M_{+}-1} 
\tilde{\Phi}[\mu]_{j}\, 
\exp \left( -\frac{(l \check{h} - \tilde{h} y_{j}/a)^2}{4\tau} \right)
\approx
\sum_{j \in J(l)} 
\tilde{\Phi}[\mu]_{j}\, 
\exp \left( -\frac{(l \check{h} - \tilde{h} y_{j}/a)^2}{4\tau} \right),
\label{eq:GaussSumApprox}
\end{align}
where $J(l)$ is the set of the indexes defined by
\begin{align}
J(l) = \{ j \mid \,  | l \check{h} - \tilde{h} y_{j}/a | \leq b \}
\label{eq:def_J_l}
\end{align}
for some $b > 0$. 
Since $y_{j}$ is defined by $\varphi$ in~\eqref{eq:DEFT_trans} as~\eqref{eq:y_j}
and monotone increasing with respect to $j$, 
the index set $J(l)$ is contained in the slightly augmented set $\tilde{J}(l)$ defined by 
\begin{align}
\tilde{J}(l) = \{ j \mid j_{\min}(l) \leq j \leq j_{\max}(l) \}, 
\label{eq:aug_J_l}
\end{align}
where
\begin{align}
j_{\min}(l) & = \max\{ j \mid l \geq \lceil (\tilde{h} y_{j}/a + b)/\check{h} \rceil \}, \label{eq:j_min}\\
j_{\max}(l) & = \max\{ j \mid l \geq \lfloor (\tilde{h} y_{j}/a - b)/\check{h} \rfloor \}. \label{eq:j_max} 
\end{align}
Figure~\ref{fig:DE_and_IndexSetJ}
shows examples of $\varphi$ and $\tilde{J}(l)$. 
As explained in Remark~\ref{rem:Comp_J_l} below,
we can obtain $\tilde{J}(l)$ for $l = -L_{-}, \ldots, L_{+}$ efficiently.
The definition~\eqref{eq:aug_J_l} of $\tilde{J}(l)$ 
means that the truncation error of approximation~\eqref{eq:GaussSumApprox} is 
$\mathrm{O}(\exp(-b^2/(4\tau)))$. 
Then, in addition to $\tau$, $a$, $\check{h}$, and $L_{\pm}$, 
we need to choose the constant $b$ 
so that both the total error 
(i.e., the sum of the discretization error for~\eqref{eq:NfftKeyRel} 
and the truncation error for~\eqref{eq:GaussSumApprox}) 
and the computational cost are reasonably small. 
For $M = M_{-} + M_{+}$ and 
a sufficiently small $\varepsilon > 0$ such as $\varepsilon = 10^{-10}$, 
we can present one possible set of their choices:%
\footnote{When $\varepsilon = 10^{-10}$, $b = 14.65\cdots$ and $\tau = 2.33\cdots$ in \eqref{eq:NFFT_para_1}.}
\begin{align}
& b = -\frac{2}{\pi} \log \varepsilon, \ 
\tau = - \frac{1}{\pi^{2}} \log \varepsilon, \ 
a = \frac{2\pi}{M}, \
\check{h} = 1, 
\label{eq:NFFT_para_1} \\
& L_{-} = b - \left\lfloor \min_{j} (\tilde{h} y_{j}/a) \right\rfloor, \
L_{+} = -L_{-} + M - 1. 
\label{eq:NFFT_para_2}
\end{align}
Under these settings, the total error of
approximations~\eqref{eq:NfftKeyRel} and~\eqref{eq:GaussSumApprox} 
is approximately $\mathrm{O}(\varepsilon)$. 
We may use larger $b$ than the value above. 
Then, we can obtain the approximate values $\tilde{\mu}_{k}$ of 
$\hat{\mu}_{k}\ (k=0,\ldots, N_{\gamma})$
as
\begin{align}
\tilde{\mu}_{k}
=
\sqrt{\frac{\pi}{\tau}}\, \e^{\tau (a k')^{2}}
\frac{\check{h}}{2\pi} 
\sum_{l = -L_{-}}^{L_{+}}
\left[
\sum_{j \in \tilde{J}(l)} 
\tilde{\Phi}[\mu]_{j}\, 
\exp \left( -\frac{(l \check{h} - \tilde{h} y_{j}/a)^2}{4\tau} \right)
\right]
\exp \left( - \i \, \frac{2\pi}{M}\, k' l \right), 
\label{eq:DEFT_final}
\end{align}
where we can use the FFT for the outer sum.
Under the settings of $a$ and $\check{h}$ in~\eqref{eq:NFFT_para_1}, 
the period of expression~\eqref{eq:DEFT_final} with respect to $k'$ is $M$. 
Therefore, we need to set 
\begin{align}
M_{-} + M_{+} = M =2N_{\gamma}
\label{eq:M}
\end{align}
to compute the accurate approximations $\tilde{\mu}_{k}\ (k=0,\ldots, N_{\gamma})$. 
Then, the computational cost is 
$\mathrm{O}(b M) + \mathrm{O}(M \log M) = \mathrm{O}(N_{\gamma} \log N_{\gamma})$. 

\begin{rem}
In~\eqref{eq:NfftKeyRel} and~\eqref{eq:DEFT_final}, 
we use the shifted index $k'$ in~\eqref{eq:k_prime} 
to avoid numerical instability in the actual computation of 
the approximate values $\tilde{\mu}_{k}$ in~\eqref{eq:DEFT_final}. 
In fact, if we use the factor $\sqrt{\pi/\tau}\, \e^{\tau (a k)^{2}}$ 
in~\eqref{eq:DEFT_final} as a usual manner of the nonuniform FFT, 
it becomes considerably large for a large $k$ and the approximation gets worse. 
\end{rem}

\begin{rem}
\label{rem:Comp_J_l}
The inequalities in~\eqref{eq:j_min} and~\eqref{eq:j_max} 
respectively defining $j_{\min}$ and $j_{\max}$ 
are nonlinear with respect to $j$, 
and it is difficult to obtain their closed forms. 
Therefore, 
noting the monotonicity of $y_{j}$ with respect to $j$, 
we use the following numerical algorithms to determine them efficiently. 
\begin{center}
\begin{minipage}[t]{0.45\linewidth}
Algorithm for $j_{\min}$\\
{\bf begin} \\
\phantom{aaaa} $j_{\min}(-L_{-} - 1) = -M_{-}$ \\
\phantom{aaaa} {\bf for } $l = -L_{-}$ to $L_{+}$\\
\phantom{aaaaaaaa} {\bf for } $j = j_{\min}(l-1)$ to $M_{+} - 1$ \\
\phantom{aaaaaaaaaaaa} {\bf if } $l \geq \lceil (\tilde{h} y_{j}/a + b)/ \check{h} \rceil$ \\
\phantom{aaaaaaaaaaaaaaaa} $j_{\min}(l) = j$ \\
\phantom{aaaaaaaaaaaa} {\bf else } \\ 
\phantom{aaaaaaaaaaaaaaaa} {\bf break } \\
\phantom{aaaaaaaaaaaa} {\bf end } \\
\phantom{aaaaaaaa} {\bf end} \\
\phantom{aaaa} {\bf end} \\
{\bf end}
\end{minipage}
\begin{minipage}[t]{0.45\linewidth}
Algorithm for $j_{\max}$\\
{\bf begin} \\
\phantom{aaaa} $j_{\max}(-L_{-} - 1) = -M_{-}$ \\
\phantom{aaaa}  {\bf for } $l = -L_{-}$ to $L_{+}$\\
\phantom{aaaaaaaa} {\bf for } $j = j_{\max}(l-1)$ to $M_{+} - 1$ \\
\phantom{aaaaaaaaaaaa} {\bf if } $l \geq \lfloor (\tilde{h} y_{j}/a - b)/\check{h} \rfloor$ \\
\phantom{aaaaaaaaaaaaaaaa} $j_{\max}(l) = j$ \\
\phantom{aaaaaaaaaaaa} {\bf else } \\ 
\phantom{aaaaaaaaaaaaaaaa} {\bf break } \\
\phantom{aaaaaaaaaaaa} {\bf end } \\
\phantom{aaaaaaaa} {\bf end} \\
\phantom{aaaa} {\bf end} \\
{\bf end}
\end{minipage}
\end{center}
\end{rem}

\subsection{Step 2: Computation of the indefinite integral of~\eqref{eq:first_FT}
in~\eqref{eq:G_alpha_1} or~\eqref{eq:G_alpha_2}}
\label{sec:SGII}

In this step, using the values
$\tilde{\mu}_{-N_{\gamma}+1}, \ldots, \tilde{\mu}_{N_{\gamma}}$, i.e., 
the approximate values of \eqref{eq:first_FT} for 
$\zeta = k \tilde{h} \ (k = -N_{\gamma}+1, \ldots, N_{\gamma})$, 
we obtain the approximate values of the indefinite integral~\eqref{eq:first_FT}
in~\eqref{eq:G_alpha_1} or~\eqref{eq:G_alpha_2} for 
$\omega = \ell \tilde{h}\ (\ell = -N+1, \ldots, N)$. 
Recall that $N$ is the integer determining the number of $x$'s 
for which we want to compute the values of the solution $p(x,t)$ in~\eqref{eq:ExprExactSol}. 
Noting~\eqref{eq:FT_conjugate}, 
we have only to compute the approximate values for 
$\omega = \ell \tilde{h} \ (\ell = 1, \ldots, N)$.

\subsubsection{Sinc-Gauss indefinite integration formula}

Let $N'$ be a positive integer. 
As a tool for computing indefinite integrals, 
we use the sinc-Gauss sampling formula~\citep{bib:KTanaka_etal_SG_2008}
for a function~$f$ on~$\mathbf{R}$
\begin{align}
f(\zeta) 
\approx 
\mathcal{T}_{N', \tilde{h}}f (\zeta)
=
\sum_{k = \lfloor \zeta/\tilde{h} \rfloor -N'+1}^{\lfloor \zeta/\tilde{h} \rfloor + N'} f(k \tilde{h})\, 
\mathop{\mathrm{sinc}} ( \zeta / \tilde{h} - k )
\exp \left( - \frac{ ( \zeta / \tilde{h} - k )^{2} }{2r^{2}} \right), 
\label{eq:SG}
\end{align}
where $\mathop{\mathrm{sinc}} ( \zeta ) = \sin(\pi \zeta) / (\pi \zeta)$. 
The error estimate of this formula is given by the following theorem,
which is a combination of the special cases of Lemmas~3.1 and~3.2 
in~\cite{bib:KTanaka_etal_SG_2008}. 

\begin{thm}[{\cite[Lemmas 3.1, 3.2]{bib:KTanaka_etal_SG_2008}}]
\label{thm:SG}
Let $f$ be an analytic and bounded function in 
$\mathcal{D}_{d} = \{ \zeta \in \mathbf{C} \mid | \mathop{\mathrm{Im}} \zeta | \leq d \}$ for some $d>0$, 
and let $\mathcal{G}_{\tilde{h}}f(\zeta) = \lim_{N' \to \infty} \mathcal{T}_{N', \tilde{h}}f(\zeta)$. 
Then, 
for a sufficiently large $N'$ and a sufficiently small $\tilde{h}>0$, 
we have the following estimates of 
the discritization error~\eqref{eq:SG_disc_estim} 
and the truncation error~\eqref{eq:SG_trun_estim} 
of approximation~\eqref{eq:SG}:
\begin{align}
& \sup_{-\infty < t < \infty} 
\left| f(t) - \mathcal{G}_{\tilde{h}}f(t) \right| 
\leq 
C r h \, \exp\left( - \frac{\pi d}{\tilde{h}} + \frac{d^{2}}{2 r^{2} \tilde{h}^{2}}\right), 
\label{eq:SG_disc_estim} \\
& \sup_{-\infty < t < \infty} 
\left| \mathcal{G}_{\tilde{h}}f(t) - \mathcal{T}_{N', \tilde{h}}f(t) \right| 
\leq 
C'\, \frac{r^{2} \e^{\frac{3}{2 r^{2}}}}{N'^{2}}\, \exp\left( - \frac{N'^{2}}{2 r^{2}} \right), 
\label{eq:SG_trun_estim} 
\end{align}
where $C$ and $C'$ are positive constants independent of $N'$, $\tilde{h}$, and $r$.
\end{thm}

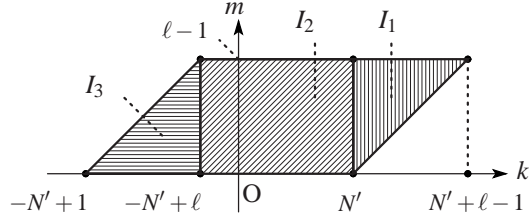
\begin{figure}[t]
\begin{center}
\unitlength 0.1in
\begin{picture}( 26.5500, 11.8000)(  0.0000,-12.0000)
%
{\color[named]{Black}{%
\special{pn 4}%
\special{sh 1}%
\special{ar 1000 1000 16 16 0  6.28318530717959E+0000}%
\special{sh 1}%
\special{ar 1000 1000 16 16 0  6.28318530717959E+0000}%
}}%
%
{\color[named]{Black}{%
\special{pn 4}%
\special{sh 1}%
\special{ar 400 1000 16 16 0  6.28318530717959E+0000}%
\special{sh 1}%
\special{ar 400 1000 16 16 0  6.28318530717959E+0000}%
}}%
%
{\color[named]{Black}{%
\special{pn 4}%
\special{sh 1}%
\special{ar 1800 1000 16 16 0  6.28318530717959E+0000}%
\special{sh 1}%
\special{ar 1800 1000 16 16 0  6.28318530717959E+0000}%
}}%
%
{\color[named]{Black}{%
\special{pn 4}%
\special{sh 1}%
\special{ar 2400 400 16 16 0  6.28318530717959E+0000}%
\special{sh 1}%
\special{ar 2400 400 16 16 0  6.28318530717959E+0000}%
}}%
%
{\color[named]{Black}{%
\special{pn 4}%
\special{sh 1}%
\special{ar 1800 400 16 16 0  6.28318530717959E+0000}%
\special{sh 1}%
\special{ar 1800 400 16 16 0  6.28318530717959E+0000}%
}}%
%
{\color[named]{Black}{%
\special{pn 4}%
\special{sh 1}%
\special{ar 1000 400 16 16 0  6.28318530717959E+0000}%
\special{sh 1}%
\special{ar 1000 400 16 16 0  6.28318530717959E+0000}%
}}%
%
{\color[named]{Black}{%
\special{pn 13}%
\special{pa 1000 400}%
\special{pa 1800 400}%
\special{pa 1800 1000}%
\special{pa 1000 1000}%
\special{pa 1000 400}%
\special{pa 1800 400}%
\special{fp}%
}}%
%
{\color[named]{Black}{%
\special{pn 4}%
\special{pa 1670 400}%
\special{pa 1200 870}%
\special{fp}%
\special{pa 1700 400}%
\special{pa 1200 900}%
\special{fp}%
\special{pa 1730 400}%
\special{pa 1200 930}%
\special{fp}%
\special{pa 1760 400}%
\special{pa 1200 960}%
\special{fp}%
\special{pa 1786 406}%
\special{pa 1200 990}%
\special{fp}%
\special{pa 1796 426}%
\special{pa 1226 996}%
\special{fp}%
\special{pa 1796 456}%
\special{pa 1256 996}%
\special{fp}%
\special{pa 1796 486}%
\special{pa 1286 996}%
\special{fp}%
\special{pa 1796 516}%
\special{pa 1316 996}%
\special{fp}%
\special{pa 1796 546}%
\special{pa 1346 996}%
\special{fp}%
\special{pa 1796 576}%
\special{pa 1376 996}%
\special{fp}%
\special{pa 1796 606}%
\special{pa 1406 996}%
\special{fp}%
\special{pa 1796 636}%
\special{pa 1436 996}%
\special{fp}%
\special{pa 1796 666}%
\special{pa 1466 996}%
\special{fp}%
\special{pa 1796 696}%
\special{pa 1496 996}%
\special{fp}%
\special{pa 1796 726}%
\special{pa 1526 996}%
\special{fp}%
\special{pa 1796 756}%
\special{pa 1556 996}%
\special{fp}%
\special{pa 1796 786}%
\special{pa 1586 996}%
\special{fp}%
\special{pa 1796 816}%
\special{pa 1616 996}%
\special{fp}%
\special{pa 1796 846}%
\special{pa 1646 996}%
\special{fp}%
\special{pa 1796 876}%
\special{pa 1676 996}%
\special{fp}%
\special{pa 1796 906}%
\special{pa 1706 996}%
\special{fp}%
\special{pa 1796 936}%
\special{pa 1736 996}%
\special{fp}%
\special{pa 1796 966}%
\special{pa 1766 996}%
\special{fp}%
\special{pa 1640 400}%
\special{pa 1200 840}%
\special{fp}%
\special{pa 1610 400}%
\special{pa 1200 810}%
\special{fp}%
\special{pa 1580 400}%
\special{pa 1200 780}%
\special{fp}%
\special{pa 1550 400}%
\special{pa 1200 750}%
\special{fp}%
\special{pa 1520 400}%
\special{pa 1200 720}%
\special{fp}%
\special{pa 1490 400}%
\special{pa 1200 690}%
\special{fp}%
\special{pa 1460 400}%
\special{pa 1200 660}%
\special{fp}%
\special{pa 1430 400}%
\special{pa 1200 630}%
\special{fp}%
\special{pa 1400 400}%
\special{pa 1200 600}%
\special{fp}%
\special{pa 1370 400}%
\special{pa 1200 570}%
\special{fp}%
\special{pa 1340 400}%
\special{pa 1200 540}%
\special{fp}%
\special{pa 1310 400}%
\special{pa 1200 510}%
\special{fp}%
\special{pa 1280 400}%
\special{pa 1200 480}%
\special{fp}%
\special{pa 1250 400}%
\special{pa 1200 450}%
\special{fp}%
\special{pa 1220 400}%
\special{pa 1200 420}%
\special{fp}%
}}%
%
{\color[named]{Black}{%
\special{pn 4}%
\special{pa 1200 660}%
\special{pa 1000 860}%
\special{fp}%
\special{pa 1200 690}%
\special{pa 1000 890}%
\special{fp}%
\special{pa 1200 720}%
\special{pa 1000 920}%
\special{fp}%
\special{pa 1200 750}%
\special{pa 1000 950}%
\special{fp}%
\special{pa 1200 780}%
\special{pa 1000 980}%
\special{fp}%
\special{pa 1200 810}%
\special{pa 1016 996}%
\special{fp}%
\special{pa 1200 840}%
\special{pa 1046 996}%
\special{fp}%
\special{pa 1200 870}%
\special{pa 1076 996}%
\special{fp}%
\special{pa 1200 900}%
\special{pa 1106 996}%
\special{fp}%
\special{pa 1200 930}%
\special{pa 1136 996}%
\special{fp}%
\special{pa 1200 960}%
\special{pa 1166 996}%
\special{fp}%
\special{pa 1200 630}%
\special{pa 1000 830}%
\special{fp}%
\special{pa 1200 600}%
\special{pa 1000 800}%
\special{fp}%
\special{pa 1200 570}%
\special{pa 1000 770}%
\special{fp}%
\special{pa 1200 540}%
\special{pa 1000 740}%
\special{fp}%
\special{pa 1200 510}%
\special{pa 1000 710}%
\special{fp}%
\special{pa 1200 480}%
\special{pa 1000 680}%
\special{fp}%
\special{pa 1200 450}%
\special{pa 1000 650}%
\special{fp}%
\special{pa 1200 420}%
\special{pa 1000 620}%
\special{fp}%
\special{pa 1190 400}%
\special{pa 1000 590}%
\special{fp}%
\special{pa 1160 400}%
\special{pa 1000 560}%
\special{fp}%
\special{pa 1130 400}%
\special{pa 1000 530}%
\special{fp}%
\special{pa 1100 400}%
\special{pa 1000 500}%
\special{fp}%
\special{pa 1070 400}%
\special{pa 1000 470}%
\special{fp}%
\special{pa 1040 400}%
\special{pa 1000 440}%
\special{fp}%
}}%
%
{\color[named]{Black}{%
\special{pn 13}%
\special{pa 400 1000}%
\special{pa 1000 1000}%
\special{pa 1000 400}%
\special{pa 400 1000}%
\special{pa 1000 1000}%
\special{fp}%
}}%
%
{\color[named]{Black}{%
\special{pn 4}%
\special{pa 980 420}%
\special{pa 1000 420}%
\special{fp}%
\special{pa 960 440}%
\special{pa 1000 440}%
\special{fp}%
\special{pa 940 460}%
\special{pa 1000 460}%
\special{fp}%
\special{pa 920 480}%
\special{pa 1000 480}%
\special{fp}%
\special{pa 900 500}%
\special{pa 1000 500}%
\special{fp}%
\special{pa 880 520}%
\special{pa 1000 520}%
\special{fp}%
\special{pa 860 540}%
\special{pa 1000 540}%
\special{fp}%
\special{pa 840 560}%
\special{pa 1000 560}%
\special{fp}%
\special{pa 820 580}%
\special{pa 1000 580}%
\special{fp}%
\special{pa 800 600}%
\special{pa 1000 600}%
\special{fp}%
\special{pa 780 620}%
\special{pa 1000 620}%
\special{fp}%
\special{pa 760 640}%
\special{pa 1000 640}%
\special{fp}%
\special{pa 740 660}%
\special{pa 1000 660}%
\special{fp}%
\special{pa 720 680}%
\special{pa 1000 680}%
\special{fp}%
\special{pa 700 700}%
\special{pa 1000 700}%
\special{fp}%
\special{pa 680 720}%
\special{pa 1000 720}%
\special{fp}%
\special{pa 660 740}%
\special{pa 1000 740}%
\special{fp}%
\special{pa 640 760}%
\special{pa 1000 760}%
\special{fp}%
\special{pa 620 780}%
\special{pa 1000 780}%
\special{fp}%
\special{pa 600 800}%
\special{pa 1000 800}%
\special{fp}%
\special{pa 580 820}%
\special{pa 1000 820}%
\special{fp}%
\special{pa 560 840}%
\special{pa 1000 840}%
\special{fp}%
\special{pa 540 860}%
\special{pa 1000 860}%
\special{fp}%
\special{pa 520 880}%
\special{pa 1000 880}%
\special{fp}%
\special{pa 500 900}%
\special{pa 1000 900}%
\special{fp}%
\special{pa 480 920}%
\special{pa 1000 920}%
\special{fp}%
\special{pa 460 940}%
\special{pa 1000 940}%
\special{fp}%
\special{pa 440 960}%
\special{pa 1000 960}%
\special{fp}%
\special{pa 420 980}%
\special{pa 1000 980}%
\special{fp}%
}}%
%
{\color[named]{Black}{%
\special{pn 8}%
\special{pa 200 1000}%
\special{pa 2600 1000}%
\special{fp}%
\special{sh 1}%
\special{pa 2600 1000}%
\special{pa 2534 980}%
\special{pa 2548 1000}%
\special{pa 2534 1020}%
\special{pa 2600 1000}%
\special{fp}%
}}%
%
{\color[named]{Black}{%
\special{pn 8}%
\special{pa 1200 1200}%
\special{pa 1200 200}%
\special{fp}%
\special{sh 1}%
\special{pa 1200 200}%
\special{pa 1180 268}%
\special{pa 1200 254}%
\special{pa 1220 268}%
\special{pa 1200 200}%
\special{fp}%
}}%
%
{\color[named]{Black}{%
\special{pn 13}%
\special{pa 1800 400}%
\special{pa 1800 1000}%
\special{pa 2400 400}%
\special{pa 1800 400}%
\special{pa 1800 1000}%
\special{fp}%
}}%
%
{\color[named]{Black}{%
\special{pn 4}%
\special{pa 1800 1000}%
\special{pa 1800 400}%
\special{fp}%
\special{pa 1820 980}%
\special{pa 1820 400}%
\special{fp}%
\special{pa 1840 960}%
\special{pa 1840 400}%
\special{fp}%
\special{pa 1860 940}%
\special{pa 1860 400}%
\special{fp}%
\special{pa 1880 920}%
\special{pa 1880 400}%
\special{fp}%
\special{pa 1900 900}%
\special{pa 1900 400}%
\special{fp}%
\special{pa 1920 880}%
\special{pa 1920 400}%
\special{fp}%
\special{pa 1940 860}%
\special{pa 1940 400}%
\special{fp}%
\special{pa 1960 840}%
\special{pa 1960 400}%
\special{fp}%
\special{pa 1980 820}%
\special{pa 1980 400}%
\special{fp}%
\special{pa 2000 800}%
\special{pa 2000 400}%
\special{fp}%
\special{pa 2020 780}%
\special{pa 2020 400}%
\special{fp}%
\special{pa 2040 760}%
\special{pa 2040 400}%
\special{fp}%
\special{pa 2060 740}%
\special{pa 2060 400}%
\special{fp}%
\special{pa 2080 720}%
\special{pa 2080 400}%
\special{fp}%
\special{pa 2100 700}%
\special{pa 2100 400}%
\special{fp}%
\special{pa 2120 680}%
\special{pa 2120 400}%
\special{fp}%
\special{pa 2140 660}%
\special{pa 2140 400}%
\special{fp}%
\special{pa 2160 640}%
\special{pa 2160 400}%
\special{fp}%
\special{pa 2180 620}%
\special{pa 2180 400}%
\special{fp}%
\special{pa 2200 600}%
\special{pa 2200 400}%
\special{fp}%
\special{pa 2220 580}%
\special{pa 2220 400}%
\special{fp}%
\special{pa 2240 560}%
\special{pa 2240 400}%
\special{fp}%
\special{pa 2260 540}%
\special{pa 2260 400}%
\special{fp}%
\special{pa 2280 520}%
\special{pa 2280 400}%
\special{fp}%
\special{pa 2300 500}%
\special{pa 2300 400}%
\special{fp}%
\special{pa 2320 480}%
\special{pa 2320 400}%
\special{fp}%
\special{pa 2340 460}%
\special{pa 2340 400}%
\special{fp}%
\special{pa 2360 440}%
\special{pa 2360 400}%
\special{fp}%
\special{pa 2380 420}%
\special{pa 2380 400}%
\special{fp}%
}}%
%
{\color[named]{Black}{%
\special{pn 13}%
\special{pa 2400 1000}%
\special{pa 2400 400}%
\special{dt 0.045}%
}}%
%
{\color[named]{Black}{%
\special{pn 4}%
\special{sh 1}%
\special{ar 2400 1000 16 16 0  6.28318530717959E+0000}%
\special{sh 1}%
\special{ar 2400 1000 16 16 0  6.28318530717959E+0000}%
}}%
\put(22.0000,-12.0000){\makebox(0,0)[lb]{{\small $N'+\ell-1$}}}%
\put(26.5500,-10.3000){\makebox(0,0)[lb]{$k$}}%
\put(17.2500,-11.9500){\makebox(0,0)[lb]{{\small $N'$}}}%
\put(12.2000,-11.5000){\makebox(0,0)[lb]{$\mathrm{O}$}}%
\put(6.0000,-12.0000){\makebox(0,0)[lb]{{\small $-N'+\ell$}}}%
\put(0.0000,-12.0000){\makebox(0,0)[lb]{{\small $-N'+1$}}}%
\put(11.3000,-1.6500){\makebox(0,0)[lb]{$m$}}%
%
{\color[named]{Black}{%
\special{pn 13}%
\special{pa 1200 400}%
\special{pa 1100 300}%
\special{dt 0.045}%
}}%
%
{\color[named]{Black}{%
\special{pn 13}%
\special{pa 800 800}%
\special{pa 600 600}%
\special{dt 0.045}%
}}%
%
{\color[named]{Black}{%
\special{pn 13}%
\special{pa 1600 600}%
\special{pa 1600 300}%
\special{dt 0.045}%
}}%
%
{\color[named]{Black}{%
\special{pn 13}%
\special{pa 2000 600}%
\special{pa 2000 300}%
\special{dt 0.045}%
}}%
\put(8.0000,-3.0000){\makebox(0,0)[lb]{{\small $\ell - 1$}}}%
\put(4.0000,-6.0000){\makebox(0,0)[lb]{$I_{3}$}}%
\put(15.0000,-2.5500){\makebox(0,0)[lb]{$I_{2}$}}%
\put(19.3500,-2.5500){\makebox(0,0)[lb]{$I_{1}$}}%
\end{picture}%
\end{center}
\caption{Partition~\eqref{eq:I_I123} of the index set $I$ in~\eqref{eq:def_IndexSetI}. 
The index sets $I_{1}$, $I_{2}$, and $I_{3}$ consist of the integer lattice points 
in the corresponding regions above, respectively. 
The boundary between $I_{1}$ and $I_{2}$, 
and the boundary between $I_{2}$ and $I_{3}$ 
belong to $I_{2}$. }
\label{fig:mkSumExc}
\end{figure}

From formula~\eqref{eq:SG}, 
we derive a formula to approximate the indefinite integral of $f$ from $0$ to $\omega$
with $\omega = \ell \tilde{h}\ (\ell = 1, \ldots, N')$.
Partitioning the integral of $f$ as
\begin{align}
\int_{0}^{\ell \tilde{h}} f(\zeta)\, \d \zeta = 
\sum_{m = 0}^{\ell-1} \int_{m \tilde{h}}^{(m+1) \tilde{h}} f(\zeta)\, \d \zeta
\label{eq:part_indefint}
\end{align}
and applying formula \eqref{eq:SG} to each term of the RHS in \eqref{eq:part_indefint}, 
we have 
\begin{align}
\int_{0}^{\ell \tilde{h}} f(\zeta)\, \d \zeta
& \approx
\sum_{m = 0}^{\ell-1} 
\sum_{k = m - N'+1}^{m + N'} f(k \tilde{h})\, 
\int_{m \tilde{h}}^{(m+1) \tilde{h}}
\mathop{\mathrm{sinc}} ( \zeta / \tilde{h} - k )
\exp \left( - \frac{ ( \zeta / \tilde{h} - k )^{2} }{2r^{2}} \right)\, \d \zeta \notag \\
& =
\sum_{m = 0}^{\ell-1} 
\sum_{k = m - N' +1}^{m + N'} f(k \tilde{h})\, 
\tilde{h}\, 
( G_{r}(m+1 - k)
-
G_{r}(m - k) ), 
\label{eq:SG_indefint}
\end{align}
where 
\begin{align}
G_{r}(\nu)
= 
\int_{0}^{\nu} 
\mathop{\mathrm{sinc}} ( \eta )
\exp \left( - \frac{ \eta^{2} }{2r^{2}} \right) \d \eta.
\label{eq:def_SGIndefInt}
\end{align}
For $\ell = 1$, straightforwardly we have 
\begin{align}
\text{(RHS of \eqref{eq:SG_indefint})}
= &
\sum_{k = - N'+1}^{N'} \tilde{h}\, f(k \tilde{h})\, G_{r}(1 - k)
-
\sum_{k = - N'+1}^{N'} \tilde{h}\, f(k \tilde{h})\, G_{r}(- k) \notag \\
= & 
\sum_{k' = - N'+1}^{N'} \tilde{h}\, f((1 - k') \tilde{h})\, G_{r}(k')
-
\sum_{k = - N'+1}^{N'} \tilde{h}\, f(k \tilde{h})\, G_{r}(- k). 
\label{eq:SG_indefint_exchanged_ell1}
\end{align}
For $\ell \geq 2$, 
setting $\tilde{f}_{k} = \tilde{h}\, f(k \tilde{h})$ and 
$d_{a}^{b}(G_{r}) = G_{r}(b) -G_{r}(a)$, 
we can rewrite the RHS of~\eqref{eq:SG_indefint} as 
\begin{align}
\text{(RHS of \eqref{eq:SG_indefint})} 
= 
\sum_{m = 0}^{\ell-1} 
\sum_{k = m - N'+1}^{m + N'} 
\tilde{f}_{k}\, d_{m-k}^{m+1-k}(G_{r}) 
= 
\sum_{(m,k) \in I} \tilde{f}_{k}\, d_{m-k}^{m+1-k}(G_{r}), 
\label{eq:SumRewrite}
\end{align} 
where $I$ is the set of indexes defined by
\begin{align}
I = \bigcup_{m = 0}^{\ell - 1} \left\{ (m,k) \mid m - N'+1 \leq k \leq m + N' \right\}. 
\label{eq:def_IndexSetI}
\end{align}
Here, we partition this index set $I$ into three disjoint parts as
\begin{align}
I = I_{1} \cup I_{2} \cup I_{3}, 
\label{eq:I_I123}
\end{align}
where 
\begin{align}
I_{1} & = \bigcup_{k=N'+1}^{N'+\ell-1} \{ (m,k) \mid k - N' \leq m \leq \ell - 1 \}, \\
I_{2} & = \bigcup_{k=-N'+\ell}^{N'} \{ (m,k) \mid 0 \leq m \leq \ell - 1 \}, \\
I_{3} & = \bigcup_{k=-N'+1}^{-N'+\ell-1} \{ (m,k) \mid 0 \leq m \leq k + N' - 1 \}. 
\end{align}
The type of this partition is illustrated by Figure~\ref{fig:mkSumExc}.
If we define $S_{i}\ (i=1,2,3)$ as
\begin{align}
S_{i} = \sum_{(m,k) \in I_{i}} \tilde{f}_{k}\, d_{m-k}^{m+1-k}(G_{r}), 
\end{align}
we have  
\begin{align}
\text{(RHS of \eqref{eq:SumRewrite})}
= S_{1} + S_{2} + S_{3}, 
\label{eq:S123}
\end{align}
and 
\begin{align}
S_{1} 
& = 
\sum_{k = N'+1}^{N' + \ell - 1} 
\sum_{m = k - N'}^{\ell-1} 
\tilde{f}_{k}\, d_{m-k}^{m+1-k}(G_{r})
=
\sum_{k = N'+1}^{N' + \ell - 1} 
\tilde{f}_{k}\, d_{-N'}^{\ell - k}(G_{r}), \label{eq:S_1} \\
S_{2}
& =
\sum_{k = - N' + \ell }^{N'} 
\sum_{m = 0}^{\ell-1} 
\tilde{f}_{k}\, d_{m-k}^{m+1-k}(G_{r})
=
\sum_{k = - N' + \ell}^{N'} 
\tilde{f}_{k}\, d_{-k}^{\ell - k}(G_{r}), \label{eq:S_2} \\
S_{3}
& = 
\sum_{k = - N' + 1}^{- N' + \ell - 1} 
\sum_{m = 0}^{k+N'-1} 
\tilde{f}_{k}\, d_{m-k}^{m+1-k}(G_{r})
=
\sum_{k = - N' + 1}^{- N' + \ell - 1} 
\tilde{f}_{k}\, d_{-k}^{N'}(G_{r}). \label{eq:S_3} 
\end{align} 
Then, combining~\eqref{eq:SumRewrite} and~\eqref{eq:S123}--\eqref{eq:S_3}, we have
\begin{align}
\text{(RHS of \eqref{eq:SG_indefint})} 
= & 
\sum_{k = - N' + \ell}^{N' + \ell - 1} \tilde{f}_{k}\, G_{r}(\ell - k)
-
\sum_{k = - N' + 1}^{N'} \tilde{f}_{k}\, G_{r}(- k) \notag \\
& 
-G_{r}(-N')
\sum_{k = N'+1}^{N' + \ell - 1} \tilde{f}_{k}\, 
+ 
G_{r}(N') 
\sum_{k = - N' + 1}^{- N' + \ell - 1 } \tilde{f}_{k}.
\label{eq:SG_indefint_exchanged_ellgeq2}
\end{align} 
Letting $k ' = \ell - k$, 
we can regard the first term of \eqref{eq:SG_indefint_exchanged_ellgeq2} as a discrete convolution
\begin{align}
\sum_{k = - N' + \ell}^{N' + \ell - 1} \tilde{h}\, f(k \tilde{h})\, G_{r}(\ell - k)
= 
\sum_{k' = - N' + 1}^{N'} \tilde{h}\, f( (\ell - k') \tilde{h} )\, G_{r}(k'). 
\label{eq:SG_convol}
\end{align}
Thus, rewriting $k'$ as $k$ in 
\eqref{eq:SG_indefint_exchanged_ell1} and 
\eqref{eq:SG_convol}, 
we finally have an indefinite integration formula: 
\begin{align}
\text{(RHS of \eqref{eq:SG_indefint})} 
= 
\sum_{k = - N' + 1}^{N'} \tilde{h}\, f((\ell - k) \tilde{h})\, G_{r}(k)
-
\sum_{k = - N' + 1}^{N'} \tilde{h}\, f(k \tilde{h})\, G_{r}(- k)
+
H_{\ell, N'}, 
\label{eq:SG_indefint_exchanged_plus}
\end{align}
where
\begin{align}
H_{\ell, N'} = 
\begin{cases}
0 & (\ell = 1), \\
\displaystyle
-
G_{r}(-N')
\sum_{k = N'+1}^{N' + \ell - 1} \tilde{h}\, f(k \tilde{h})\, 
+ 
G_{r}(N') 
\sum_{k = - N' + 1}^{- N' + \ell - 1 } \tilde{h}\, f(k \tilde{h}) 
& (\ell = 2, \ldots, N').
\label{eq:SG_indefint_exchanged_plus_where}
\end{cases}
\end{align}
To obtain the values of~\eqref{eq:SG_indefint_exchanged_plus}
for $\ell = 1,\ldots, N'$, 
we need the $3N'$ values of $f(\zeta)$ for 
$\zeta = \ell \tilde{h}\ (\ell = -N',\ldots, 2N'-1)$. 

Therefore,
to compute the indefinite integral of~\eqref{eq:first_FT} in~\eqref{eq:G_alpha_1} 
for $\omega = \ell \tilde{h}\ (\ell = 1,\ldots ,N)$, 
we set 
\begin{align}
& N_{1} = 2N \label{eq: N_alpha_1} 
\end{align}
and use formula~\eqref{eq:SG_indefint_exchanged_plus}
with $N' = N$ and $f(\ell \tilde{h})$ replaced by $\tilde{\mu}_{\ell}$ for $\ell = -N,\ldots, 2N-1$. 
Furthermore, for the integral of~\eqref{eq:first_FT} in~\eqref{eq:G_alpha_2}, we set
\begin{align}
& N_{2} = 4N \label{eq: N_alpha_2}
\end{align}
and use formula~\eqref{eq:SG_indefint_exchanged_plus} twice 
with $N' = 2N$ for the first time and $N' = N$ for the second time. 

In terms of computational time, 
note that the second term in~\eqref{eq:SG_indefint_exchanged_plus} 
and $H_{\ell, N'}\ (\ell = 2,\ldots, N')$ in~\eqref{eq:SG_indefint_exchanged_plus_where} 
can be computed in $\mathrm{O}(N)$ time when $N'=N$ or $N'=2N$. 
Then, what remains is to speed up 
the computation of the discrete convolution of the first term in \eqref{eq:SG_indefint_exchanged_plus}. 
Extending the sum to a convolution with length $4N'$ and 
using the FFT as shown in Section~\ref{sec:convol_FFT} below, 
we can compute the discrete convolution of the first term 
in $\mathrm{O}(N \log N)$ time when $N'=N$ or $N'=2N$.

\subsubsection{Fast computation of the convolution using the FFT}
\label{sec:convol_FFT}

Consider the first term in \eqref{eq:SG_indefint_exchanged_plus} 
with $f((\ell - k) \tilde{h})$ replaced by $\tilde{\mu}_{\ell - k}$: 
\begin{align}
\tilde{h} \sum_{k = - N' + 1}^{N'} \tilde{\mu}_{\ell - k}\, G_{r}(k). 
\end{align}
We compute this convolution for $\ell = -N'+1, \ldots, N'$, 
although its values for $\ell = -N', \ldots, 0$ are not required. 
To use the FFT for this computation, 
we define the sequence $\{ g[r]_{k} \}_{k=-2N'+1}^{2N'}$ as 
\begin{align}
g[r]_{k} = 
\begin{cases}
G_{r}(k) & (-N'+1 \leq k \leq N'), \\
0 & (-2N'+1\leq k \leq -N',\ N'+1 \leq k \leq 2N').
\end{cases}
\label{eq:def_g_r}
\end{align}
Then, we have
\begin{align}
\tilde{h} \sum_{k = - N' + 1}^{N'} \tilde{\mu}_{\ell - k}\, G_{r}(k)
=
\frac{\tilde{h} }{4N'}
\sum_{m = -2N'+1}^{2N'} \mathrm{DFT}[\tilde{\mu}]_{m}\, \mathrm{DFT}[g[r]]_{m}\,  
\e^{\i\, \frac{2\pi}{4N'}\, \ell\, m} 
\label{eq:DFT_IDFT}
\end{align}
for $\ell = -N'+1, \ldots, N'$, where 
\begin{align}
\mathrm{DFT}[\tilde{\mu}]_{m}
& =
\sum_{k_{1} = -2N'+1}^{2N'} \tilde{\mu}_{k_{1}}\, \e^{-\i\, \frac{2\pi}{4N'}\, k_{1}\, m}
\quad (m=-2N'+1,\ldots, 2N'), \\
\mathrm{DFT}[g[r]]_{m}
& =
\sum_{k_{2} = -2N'+1}^{2N'} g[r]_{k_{2}}\, \e^{-\i\, \frac{2\pi}{4N'}\, k_{2}\, m} 
\quad (m=-2N'+1,\ldots, 2N').
\end{align}
For the computation of~\eqref{eq:DFT_IDFT}, 
we need the values of $G_{r}(k)$ in~\eqref{eq:def_SGIndefInt}. 
In fact, they can also be computed accurately 
by a Fourier-based method and the FFT 
as presented in Appendix A.
Therefore, we can compute~\eqref{eq:DFT_IDFT} by the FFT
in $\mathrm{O}(N \log N)$ time when $N'=N$ or $N'=2N$. 

\begin{rem}
\label{rem:ErrSGindef}
Since the indefinite integration formula~\eqref{eq:SG_indefint_exchanged_plus} 
is derived from the sinc-Gauss sampling formula~\eqref{eq:SG}, 
formula~\eqref{eq:SG_indefint_exchanged_plus} 
inherits the error of formula~\eqref{eq:SG} estimated in Theorem~\ref{thm:SG}. 
In particular, the error of formula~\eqref{eq:SG_indefint_exchanged_plus}
is bounded by one of formula~\eqref{eq:SG} multiplied by $N' \tilde{h}$. 
According to Theorem~\ref{thm:SG}, 
the optimal settings of $\tilde{h}$ and $r$ for fixed $N'$
are $\tilde{h} = d/N'$ and $r = \sqrt{N'/\pi}$, respectively, 
and the total error of formula~\eqref{eq:SG} under these settings is 
$\mathrm{O}(\sqrt{1/N'} \exp(-(\pi/2) N'))$~\citep[Theorem 3.3]{bib:KTanaka_etal_SG_2008}. 
In this paper, however, 
we use the grid spacing $\tilde{h}$ determined by~\eqref{eq:EulerFTpara} 
in Theorem~\ref{thm:EulerFTError} in Section~\ref{sec:EulerFT} below, 
which results in $\tilde{h} = \mathrm{O}(\sqrt{1/N'})$. 
This choice gives priority to the theoretical settings of the parameters 
in formula~\eqref{eq:LastEulerFFT} 
for the inverse Fourier transform~\eqref{eq:ExprExactSol} in Step~3. 
Then, this $\tilde{h}$ and $r = \sqrt{N'/\pi}$ 
yield the total error $\mathrm{O}(\exp(-c \sqrt{N'}))$ of formula~\eqref{eq:SG}, 
which has a similar exponential part 
to error~\eqref{eq:EulerFTerror} of formula~\eqref{eq:LastEulerFFT} 
with respect to $N$ when $N' = N$ or $N' = 2N$. 
\end{rem}

\begin{rem}
The new formula~\eqref{eq:SG_indefint_exchanged_plus} 
is introduced for fast computation of highly accurate approximations 
of an indefinite integral on the equispaced grid 
from the values of the integrand on the same grid. 
Among the traditional quadrature formulas, 
the Newton-Cotes formulas enable such computation. 
However, these formulas have errors 
$\mathrm{O}(\tilde{h}^{\kappa})$ for some $\kappa > 0$,  
and become algebraic with respect to $N'$ 
when $\tilde{h} = \mathrm{O}((N')^{-\lambda})$ for some $\lambda > 0$, 
whereas the formula~\eqref{eq:SG_indefint_exchanged_plus} realizes the exponential convergence 
as shown in Remark~\ref{rem:ErrSGindef}. 
\end{rem}

\begin{rem}
The partition~\eqref{eq:I_I123} of the index set $I$ in~\eqref{eq:def_IndexSetI}
shown by Figure~\ref{fig:mkSumExc}
is the key to the derivation of formula~\eqref{eq:SG_indefint_exchanged_plus}. 
A similar but different idea is proposed in~\cite{bib:HaleTownsend}
for the computation of the convolution of functions. 
\end{rem}

\subsection{Step 3: Computation of the inverse Fourier transform}
\label{sec:EulerFT}

Let $[\tilde{\mathcal{G}}_{\gamma} \mu]_{N}(\ell \tilde{h})
\ (\ell = -N+1, \ldots, N)$ denote 
the approximations of~\eqref{eq:G_alpha_1} or~\eqref{eq:G_alpha_2}
computed in Step~2. 
In order to approximate the inverse Fourier transform~\eqref{eq:ExprExactSol}, 
we use the formula for the Fourier transform with a continuous Euler transform 
introduced by \cite{bib:OouraEuler2001}. 
Define $w(y; p, q)$ by
\begin{align}
w(\xi; p, q) = \frac{1}{2} \mathop{\mathrm{erfc}} \left( \frac{\xi}{p} - q \right), 
\label{eq:def_w}
\end{align}
where $\mathop{\mathrm{erfc}}$ is the complementary error function defined as
\begin{align}
\mathop{\mathrm{erfc}}(\xi) = \frac{2}{\sqrt{\pi}} \int_{\xi}^{\infty} \exp(-t^{2})\, \mathrm{d}t.
\end{align}
Using $w(\xi; p, q)$ as a weight function, 
we consider the following approximations of~\eqref{eq:ExprExactSol}: 
\begin{align}
\frac{1}{2\pi} \int_{-\infty}^{\infty} 
\exp( t\, [\mathcal{G}_{\gamma} \mu](\omega) ) \, 
\e^{\i\, x\, \omega}\, 
\d \omega 
& \approx
\frac{1}{2\pi} \int_{-\infty}^{\infty} 
w(| \omega |; p,q)
\exp( t\, [\mathcal{G}_{\gamma} \mu](\omega) ) \, 
\e^{\i\, x\, \omega}\, 
\d \omega \notag \\
& \approx 
\frac{\tilde{h}}{2\pi}\,
\sum_{\ell = -N+1}^{N}
w(| \ell \tilde{h} |; p,q) 
\exp( t\, [\tilde{\mathcal{G}}_{\gamma} \mu]_{N}(\ell \tilde{h}) ) \, 
\e^{\i\, x\, \ell \tilde{h}}.
\label{eq:LastEulerFFT}
\end{align}
The formula~\eqref{eq:LastEulerFFT} is 
the formula for the Fourier transform with a continuous Euler transform. 
The role of the function $w(|\omega|; p, q)$ 
is to realize the rapid decay of the integrand as $|\omega| \to \infty$ on $\mathbf{R}$. 
Then, we need to compute the values of~\eqref{eq:LastEulerFFT} 
for $x = n \hat{h}\ (n = -N+1,\ldots, N)$. 
Substituting this expression of $x$ into the factor $\e^{\i\, x\, \ell \tilde{h}}$ 
in~\eqref{eq:LastEulerFFT}, we have
\begin{align}
\exp( \i\, x\, \ell \tilde{h}) = \exp( \i\, \tilde{h} \hat{h}\, \ell n).
\end{align}
Unless the Nyquist condition $\tilde{h} \hat{h} = \pi/N$ holds, 
we cannot apply the FFT directly to the computation of~\eqref{eq:LastEulerFFT}. 
Therefore, for this computation, 
we use the fractional FFT~\citep{bib:BailSwarFRFT1991} 
that enables fast computation of the DFT without the Nyquist condition. 
Then, we can compute~\eqref{eq:LastEulerFFT} in $\mathrm{O}(N \log N)$ time. 
The details of the fractional FFT combined with 
the formula~\eqref{eq:LastEulerFFT} is explained in~\cite{bib:KTanaka_EulerFFT_2014}.
The error bound of the formula~\eqref{eq:LastEulerFFT} 
is given by Theorem~4 in~\cite{bib:KTanaka_EulerFFT_2014}. 

\begin{thm}[{\cite[Theorem~4]{bib:KTanaka_EulerFFT_2014}}]
\label{thm:EulerFTError}
Let $f$ be a function analytic and bounded in $\mathcal{T}'_{\theta} \cup \mathcal{D}_{d}$ 
for some $\theta$ with $0< \theta <\pi/2$ and $d>0$, 
where
\begin{align}
\mathcal{T}'_{\theta} 
& = 
\{ z\in \mathbf{C} \mid 
| \mathop{\mathrm{arg}} z | < \theta \text{ or } 
| \pi - \mathop{\mathrm{arg}} z | < \theta \}, \label{eq:domT} \\
\mathcal{D}_{d}
& = 
\{ z\in \mathbf{C} \mid 
| \mathop{\mathrm{Im}} z | < d \}
\label{eq:domD}. 
\end{align}
Moreover, assume that 
\begin{align}
\lim_{R \to \infty} \max_{-\theta \leq \phi \leq \theta} | f(\pm R + \i\, (\tan \phi) R ) | = 0
\end{align}
and $f$ is square integrable on $\mathbf{R}$. 
Let $x_{l}$ and $x_{u}$ be real numbers with $0 < x_{l} < x_{u}$ and $x_{l}/x_{u} \leq \min\{ \tan \theta, 1/2\}$. 
Then, for any $x$ with $x_{l} \leq |x| \leq x_{u}$ and a sufficiently large integer $N>0$, 
defining $h$, $p$, $q$ by
\begin{align}
\tilde{h} = \sqrt{\frac{2\pi d\, (x_{l} + x_{u})}{x_{l}^{2} N}}, \quad 
p = \sqrt{\frac{N \tilde{h}}{x_{l}}}, \quad
q = \sqrt{\frac{x_{l} N \tilde{h}}{4}}, 
\label{eq:EulerFTpara}
\end{align}
we have 
\begin{align}
\left|
\int_{-\infty}^{\infty} f(\omega)\, \e^{\i\, x\, \omega}\, \d \omega
- 
\tilde{h} \sum_{\ell = -N+1}^{N}
w(| \ell \tilde{h} |; p,q) 
f(\ell \tilde{h})\, \e^{\i\, x\, \ell \tilde{h}}
\right|
= 
\mathrm{O}
\left[
\sqrt{N}\, \exp \left( - \sqrt{\frac{\pi d\, x_{l}^{2} N}{2(x_{l} + x_{u})}} \right)
\right].
\label{eq:EulerFTerror}
\end{align}
\end{thm}

\begin{rem}
\label{rem:x_l}
Since the assumption of Theorem~\ref{thm:EulerFTError} 
includes the case that $f$ is not absolute integrable on~$\mathbf{R}$, 
the Fourier transform of $f$ may be discontinuous or non-smooth at the origin $x=0$.
Then, we consider the positive lower bound $x_{l}$ of the absolute value of $x$ 
to avoid the error estimate around the origin. 
See \cite{bib:KTanaka_EulerFFT_2014} for the details of the error estimate. 
\end{rem}

\begin{rem}
We use the different formulas~\eqref{eq:DE_FT} and~\eqref{eq:LastEulerFFT}
for the Fourier transforms~\eqref{eq:first_FT} and~\eqref{eq:ExprExactSol}, 
respectively. 
This is because the former Fourier transform~\eqref{eq:first_FT} is the integral on 
the semi-infinite interval $[0, \infty)$ 
whereas the latter Fourier transform~\eqref{eq:ExprExactSol} 
is one on the infinite interval $(-\infty, \infty)$. 
Further, we can also use formula~\eqref{eq:DE_FT} for~\eqref{eq:ExprExactSol} 
by partitioning the interval $(-\infty, \infty)$ to $(-\infty, 0]$ and $[0, \infty)$. 
However, we give priority to 
a brief implementation and lower computational cost 
of formula~\eqref{eq:LastEulerFFT}. 
\end{rem}

\section{Numerical examples}
\label{sec:NumEx}

In this section, 
we apply the proposed method 
to two examples of the PIDE \eqref{eq:PIDE_LM}
with initial condition $p(x,0) = \delta(x)$. 
Let $K_{v}$ be the modified Bessel function of the second kind. 

\begin{ex}[Variance gamma (VG) process~\citep{bib:Apple_LevyText_2009}]
\label{ex:VG}
Setting $\gamma = 1$ and $\mu(y) = \e^{- y}$ in~\eqref{eq:Lv_measure}, 
we have the L\'evy measure
\begin{align}
\nu(\d y) = \frac{1}{|y|}\, e^{-|y|}\, \d y. 
\label{eq:sym_VG_meas}
\end{align}
The L\'evy process described by this measure is 
the symmetric VG process. 
The exact solution of \eqref{eq:PIDE_LM} with measure~\eqref{eq:sym_VG_meas} 
is written in the form
\begin{align}
p(x, t) = 
\left( \frac{|x|}{2} \right)^{t-1/2} \frac{K_{1/2-t} (|x|)}{\sqrt{\pi}\, \Gamma(t)}. 
\label{eq:exact_sol_VG}
\end{align}
\end{ex}

\begin{ex}[Normal inverse Gaussian (NIG) process~\citep{bib:Apple_LevyText_2009}]
\label{ex:NIG}
Setting $\gamma = 2$ and $\mu(y) = y\, K_{1}(y)/\pi$ in~\eqref{eq:Lv_measure}, 
we have the L\'evy measure
\begin{align}
\nu(\d y) = \frac{1}{\pi |y|}\, K_{1}(|y|)\, \d y. 
\label{eq:sym_NIG_meas}
\end{align}
The L\'evy process described by this measure is 
the symmetric NIG process. 
The exact solution of \eqref{eq:PIDE_LM} with measure~\eqref{eq:sym_NIG_meas} 
is written in the form
\begin{align}
p(x, t) = t\, \e^{t}\, \frac{K_{1}( \sqrt{x^2 + t^2} )}{ \pi \sqrt{x^2 + t^2} }. 
\label{eq:exact_sol_NIG}
\end{align}
\end{ex}

Using the proposed method, 
we compute the numerical solutions of these examples
for $x \in [-5,5]$ and $t = 1,2,3$. 
Then, $x_{u}$ in Theorem~\ref{thm:EulerFTError} should be set as $x_{u} = 5$. 
In addition, we choose $x_{l} = 2$. 
To set equispaced grids on $[-5,5]$,
we consider the sampling points $x = n \hat{h}\ (n=-N+1,\ldots, N)$ with $\hat{h} = 5/N$
for 
\begin{align}
N = 2^{i-i_{\gamma}}\quad (i=7,\ldots, 12),
\label{eq:set_N}
\end{align}
where $i_{1} = 2$ for Example~\ref{ex:VG} and $i_{2} = 3$ for Example~\ref{ex:NIG}. 
The other parameters in the proposed method are determined as described below.
\begin{description}
\setlength{\parskip}{0pt}
\setlength{\itemsep}{3pt}
\item[Step 1] 
For Example~\ref{ex:VG}, 
$M = 4N$ and $N_{1} = 2N$ according to~\eqref{eq:M} and~\eqref{eq: N_alpha_1}. 
For Example~\ref{ex:NIG}, 
$M = 8N$ and $N_{2} = 4N$ according to~\eqref{eq:M} and~\eqref{eq: N_alpha_2}. 
In addition, $\tilde{h} = \sqrt{7\pi/(2N)}$ according to Theorem~\ref{thm:EulerFTError}, 
where we use $d=1$.
The other parameters required in Step~1 are as follows:
\begin{align}
& M_{-} = M/2,\ M_{+} = M-M_{-},\ h = \log(10^3 M)/M, 
\label{eq:para_DEFT_1} \\ 
& \zeta_{0} = [\text{the values of \eqref{eq:set_zeta0_1} and~\eqref{eq:set_zeta0_2}}],\ 
\beta = 0.25,\ 
\alpha = [\text{the value of \eqref{eq:def_DEFT_alpha}}],\ 
\label{eq:para_DEFT_2} \\
& b = 20,\ 
(\tau, a, \check{h}, L_{\pm}) =  
[\text{the set of the values of \eqref{eq:NFFT_para_1} and \eqref{eq:NFFT_para_2} for $\varepsilon = 10^{-10}$}]. 
\label{eq:para_NFFT} 
\end{align}
\item[Step 2] 
According to Theorem~\ref{thm:SG} and Remark~\ref{rem:ErrSGindef}, 
we set $r = \sqrt{N'/\pi}$ in the Gaussian kernel in~\eqref{eq:def_SGIndefInt}, 
where $N'=N$ for Example~1, and $N'=2N$ and $N'=N$ 
in the first and second application of the indefinite integral formula, respectively, 
for Example~2.
\item[Step 3] 
According to Theorem~\ref{thm:EulerFTError}, 
we set $p = q = \sqrt{N \tilde{h}/2}$ because $x_{l} = 2$. 
\end{description}
In~\eqref{eq:para_DEFT_1}, 
$h$ is not set based on a theoretical criterion 
but it is determined experimentally in reference to the settings 
in the DE formulas for definite integration~\citep{bib:KTanaka_etal_DE_2009}. 
All computations are performed through MATLAB R2013a programs 
with double precision floating point arithmetic on a PC with a
$3.0$ GHz CPU and $2.0$ GB RAM. 
The Matlab codes used for the computations are exposed 
on web page~\cite{bib:KTanaka_Matlab_2014}. 

\begin{figure}[H]
\begin{center}

\begin{minipage}[t]{0.45\linewidth}
\includegraphics[width = \linewidth]{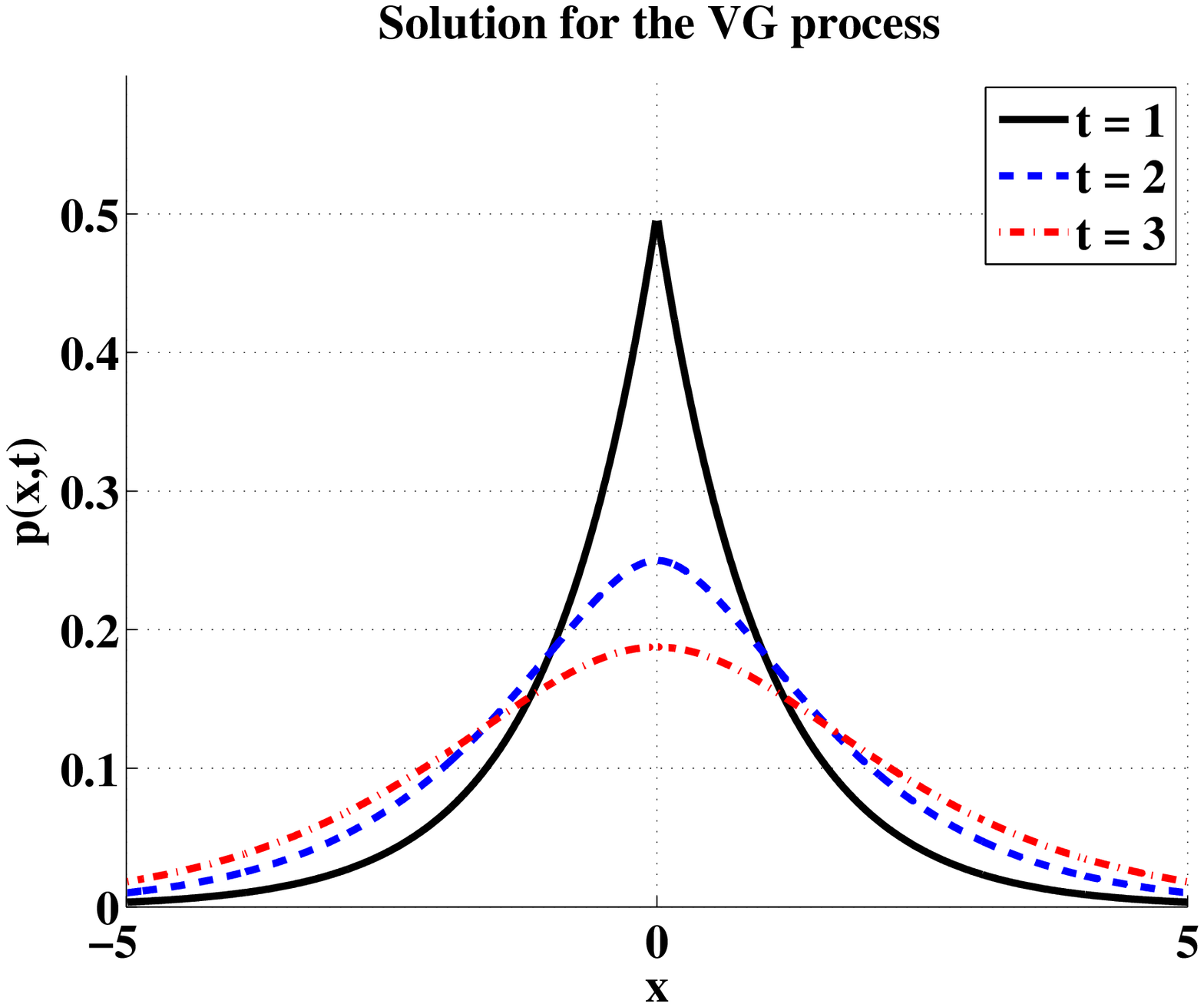}
\caption{Exact solutions $p(x,t)$ in~\eqref{eq:exact_sol_VG} of Example~\ref{ex:VG}}
\label{fig:sol_VG}
\end{minipage}\quad 
\begin{minipage}[t]{0.45\linewidth}
\includegraphics[width = \linewidth]{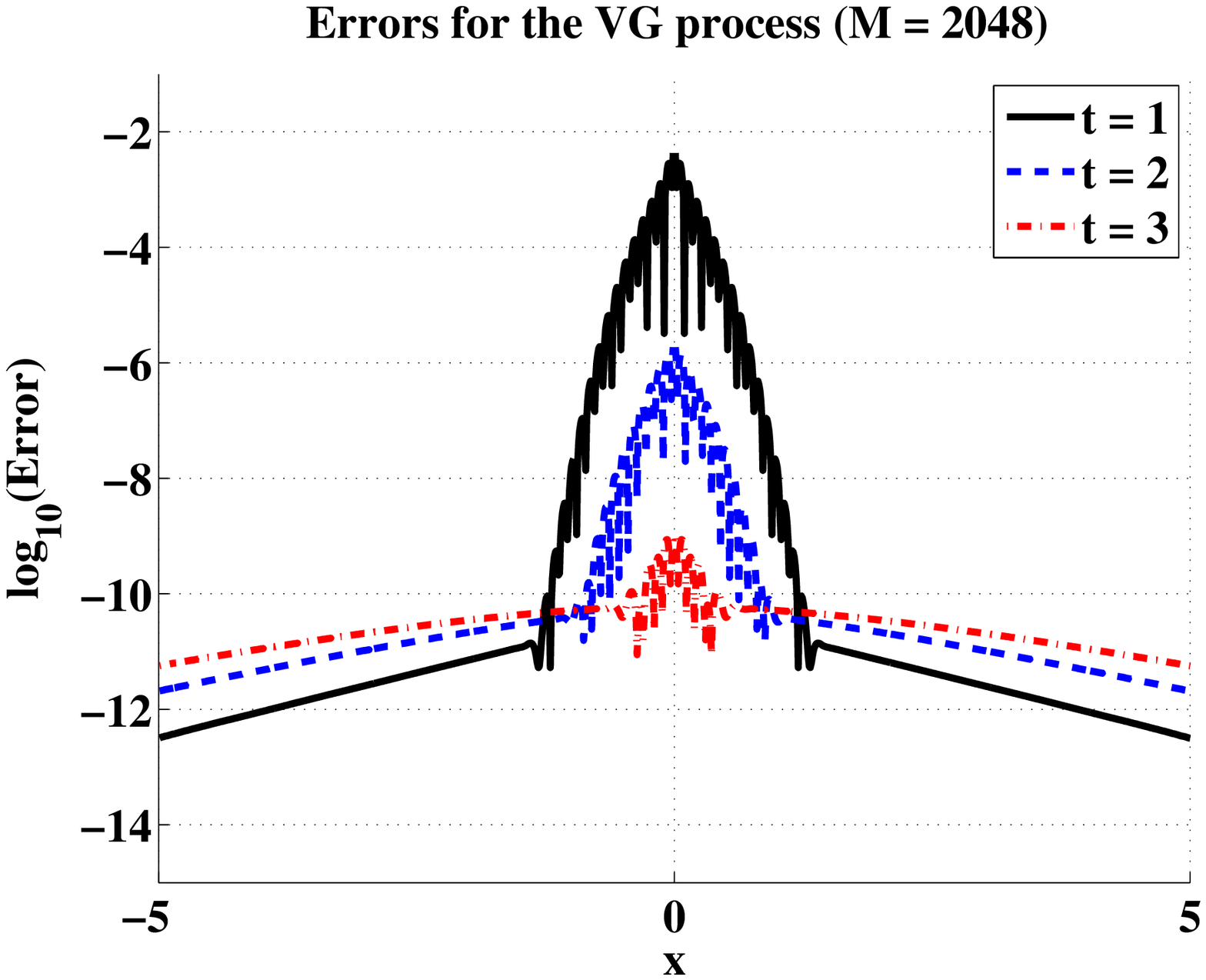}
\caption{Errors of the numerical solutions for Example~\ref{ex:VG} 
when $t = 1,2,3$ and $M=2^{11}$ (i.e., $N=2^{9}$)}
\label{fig:err_VG_M_2_11}
\end{minipage}

\end{center}
\end{figure}

\begin{figure}[H]
\begin{center}

\begin{minipage}[t]{0.45\linewidth}
\includegraphics[width = \linewidth]{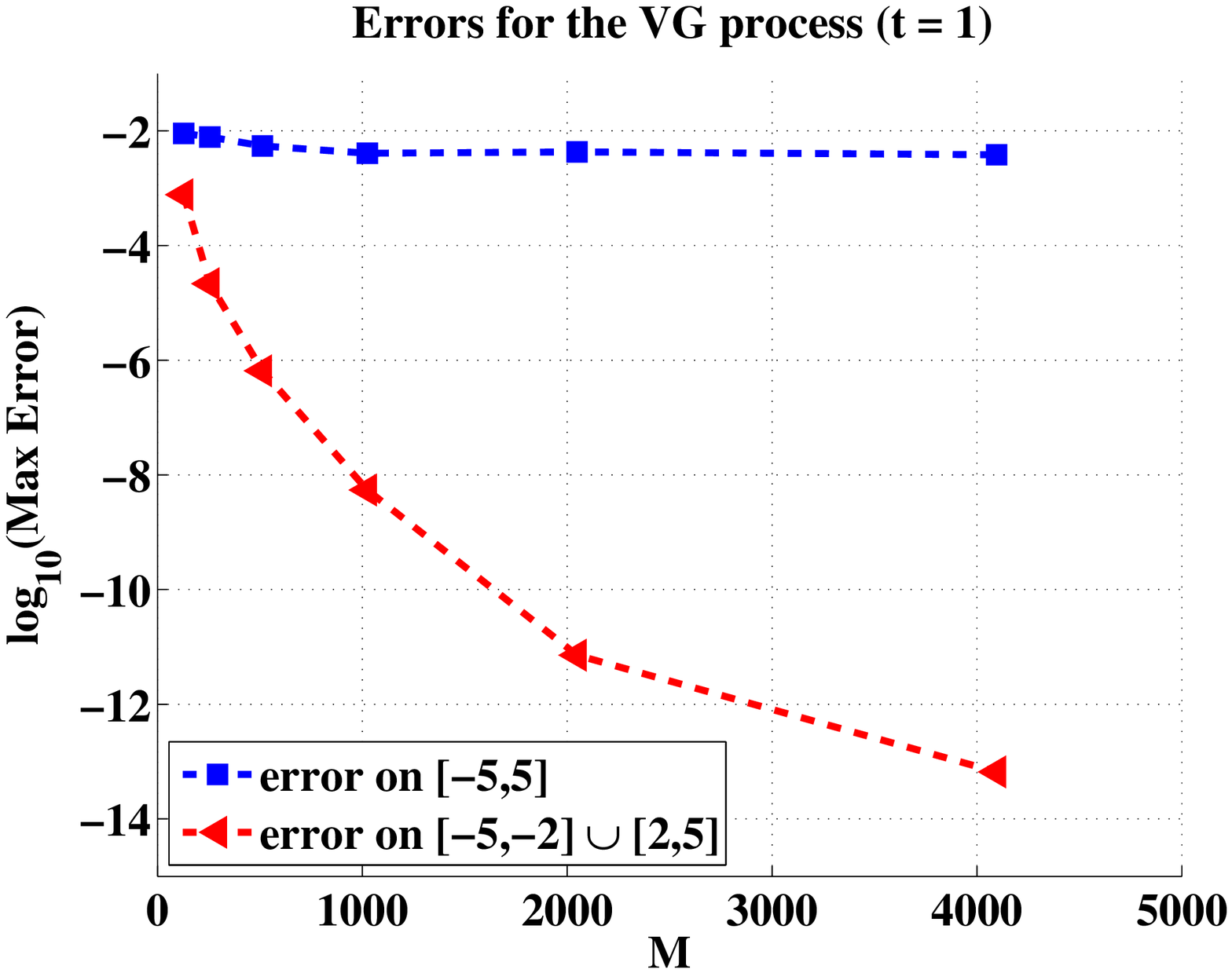}
\caption{Errors of the numerical solutions 
on $[-5, 5]$ and $[-5, -2] \cup [2, 5]$ for Example~\ref{ex:VG} 
when $t = 1$ and $M = 2^{7}, \ldots,  2^{12}$ (i.e., $N = 2^{5}, \ldots,  2^{10}$)}
\label{fig:err_VG_t1}
\end{minipage}\quad 
\begin{minipage}[t]{0.45\linewidth}
\includegraphics[width = \linewidth]{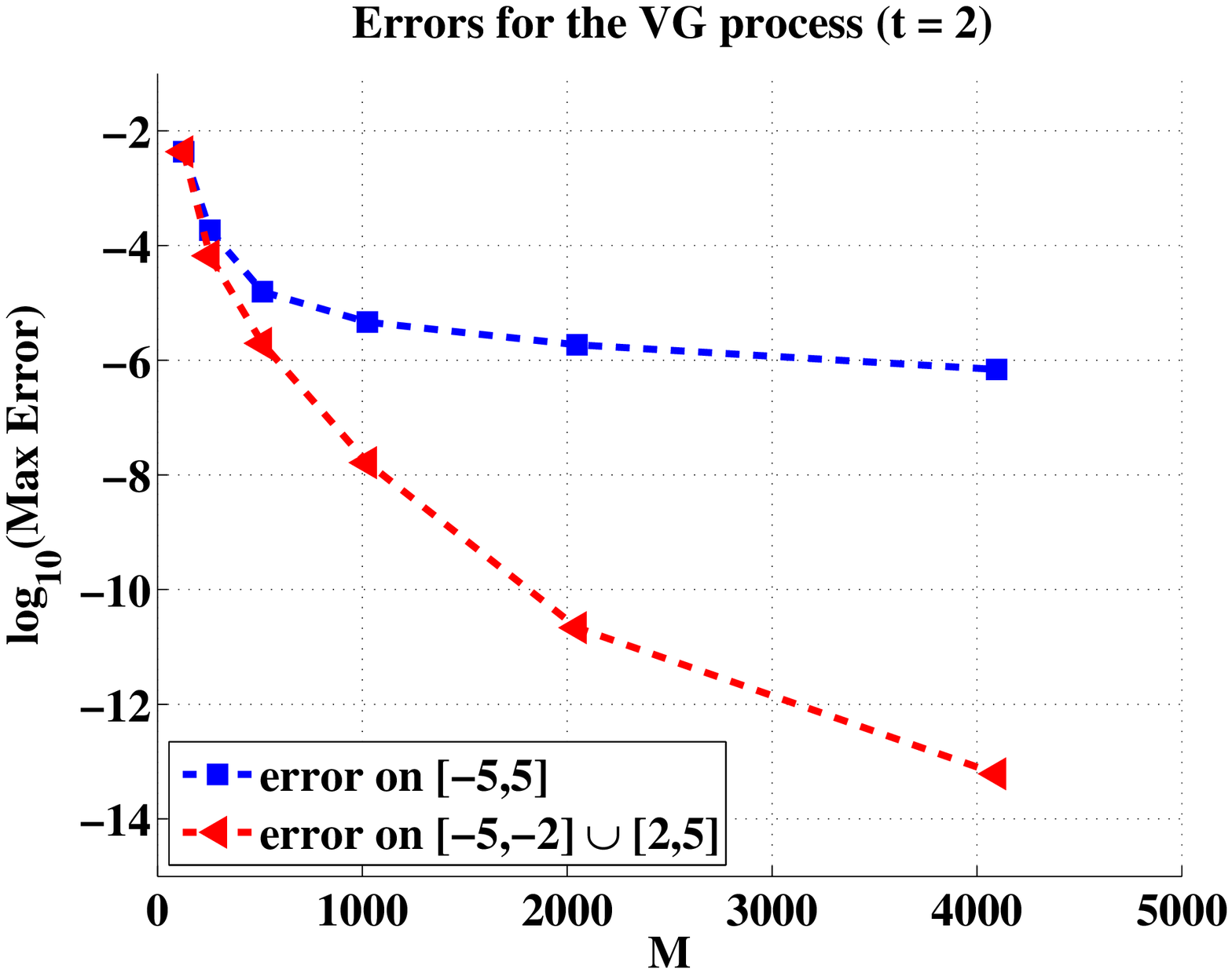}
\caption{Errors of the numerical solutions 
on $[-5, 5]$ and $[-5, -2] \cup [2, 5]$ for Example~\ref{ex:VG} 
when $t = 2$ and $M = 2^{7}, \ldots,  2^{12}$ (i.e., $N = 2^{5}, \ldots,  2^{10}$)}
\label{fig:err_VG_t2}
\end{minipage}

\end{center}
\end{figure}

\begin{figure}[H]
\begin{center}

\begin{minipage}[t]{0.45\linewidth}
\includegraphics[width = \linewidth]{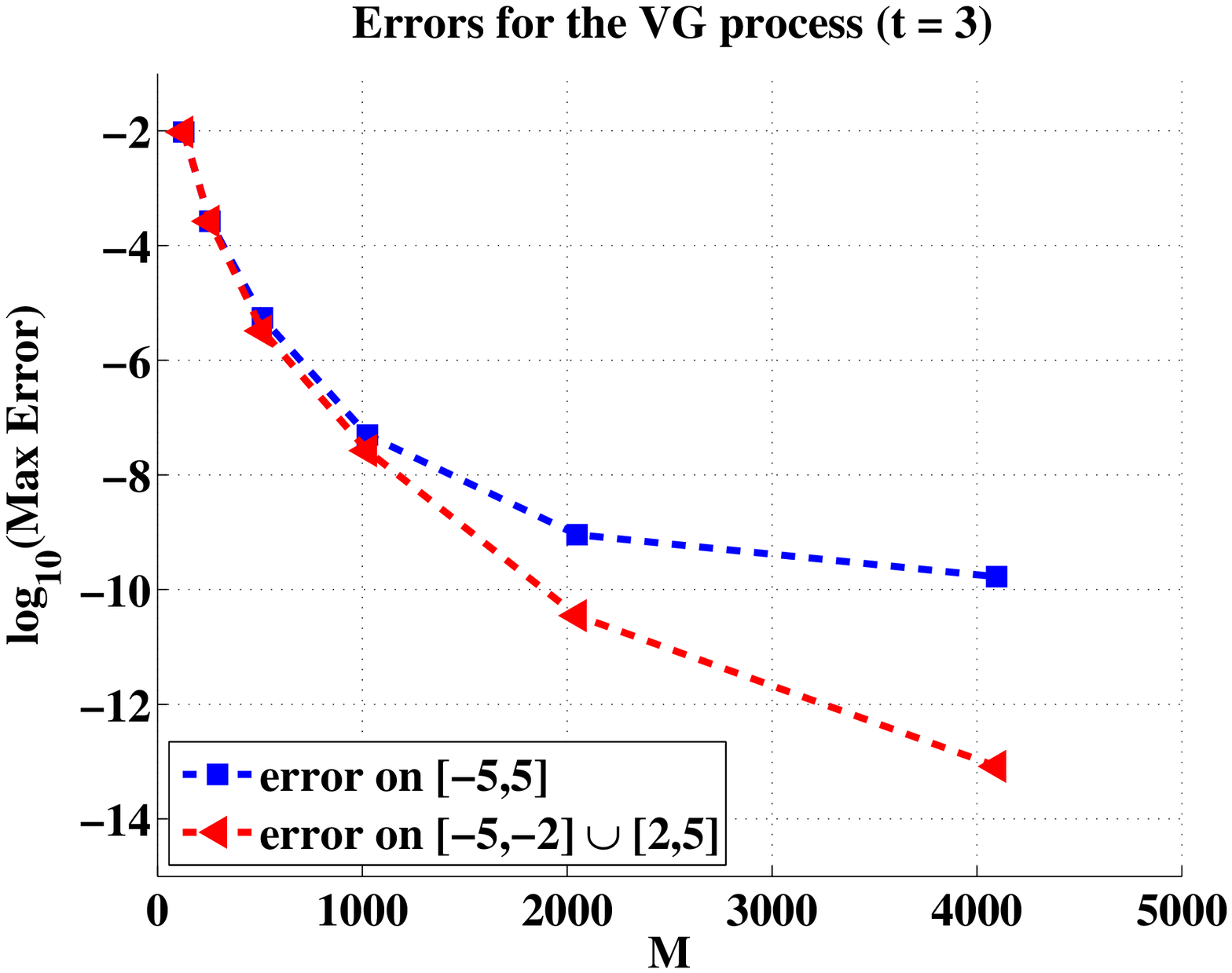}
\caption{Errors of the numerical solutions 
on $[-5, 5]$ and $[-5, -2] \cup [2, 5]$ for Example~\ref{ex:VG} 
when $t = 3$ and $M = 2^{7}, \ldots,  2^{12}$ (i.e., $N = 2^{5}, \ldots,  2^{10}$)}
\label{fig:err_VG_t3}
\end{minipage}\quad 
\begin{minipage}[t]{0.45\linewidth}
\includegraphics[width = \linewidth]{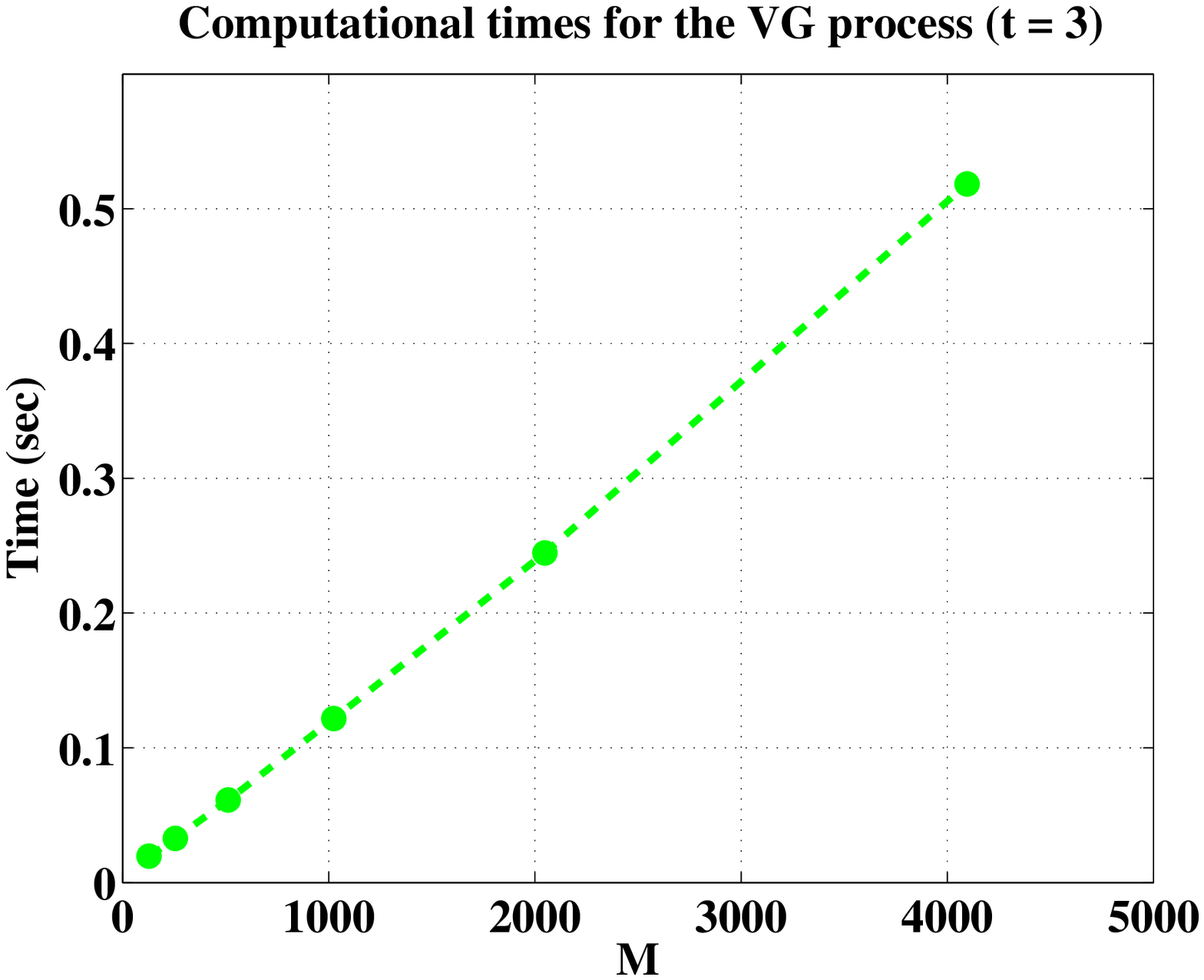}
\caption{Computational times for Example~\ref{ex:VG} when $t = 3$
and $M = 2^{7}, \ldots,  2^{12}$ (i.e., $N = 2^{5}, \ldots,  2^{10}$)}
\label{fig:time_VG_t3}
\end{minipage}

\end{center}
\end{figure}

\begin{figure}[H]
\begin{center}

\begin{minipage}[t]{0.45\linewidth}
\includegraphics[width = \linewidth]{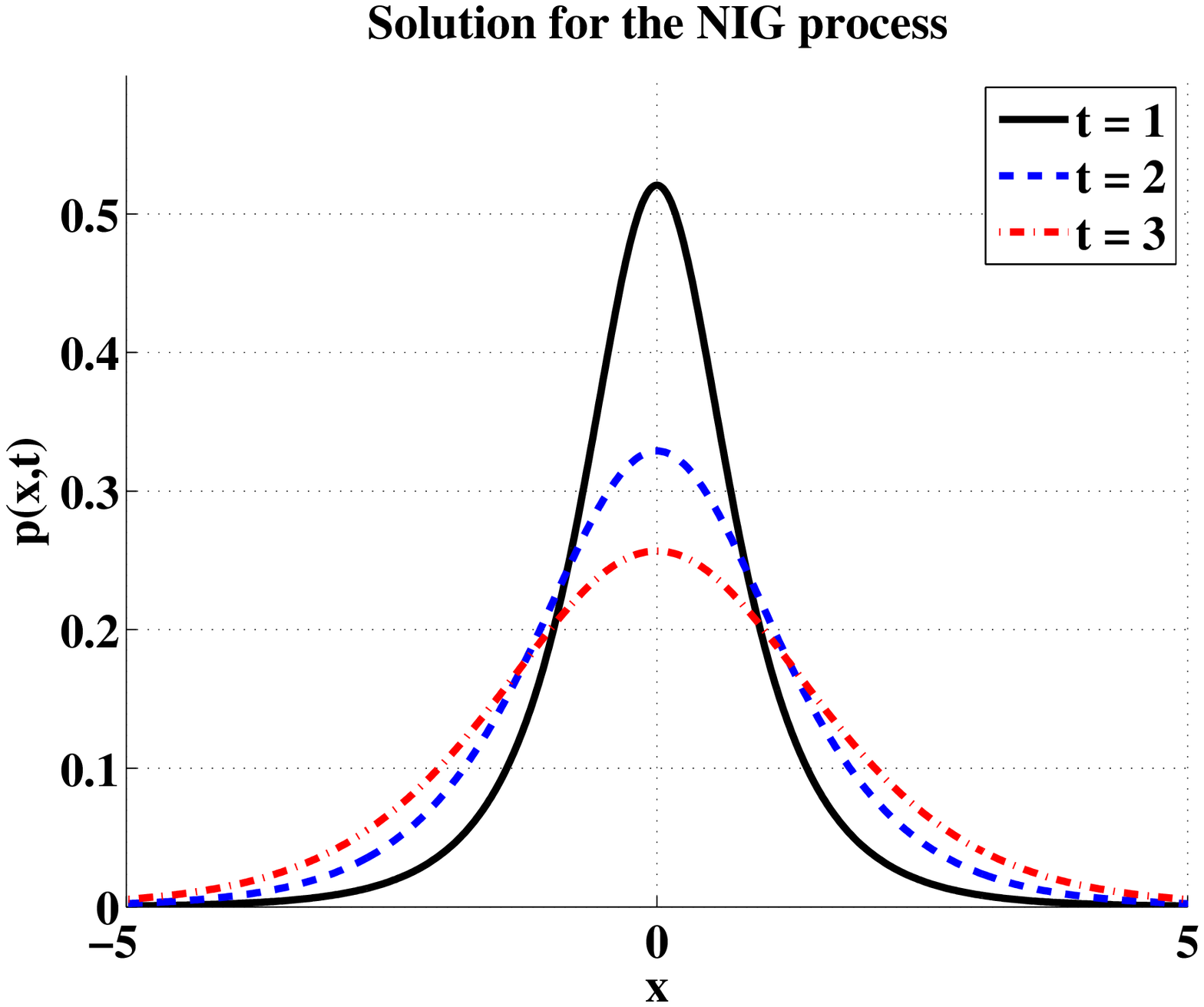}
\caption{Exact solutions $p(x,t)$ in~\eqref{eq:exact_sol_NIG} of Example~\ref{ex:NIG}}
\label{fig:sol_NIG}
\end{minipage}\quad
\begin{minipage}[t]{0.45\linewidth}
\includegraphics[width = \linewidth]{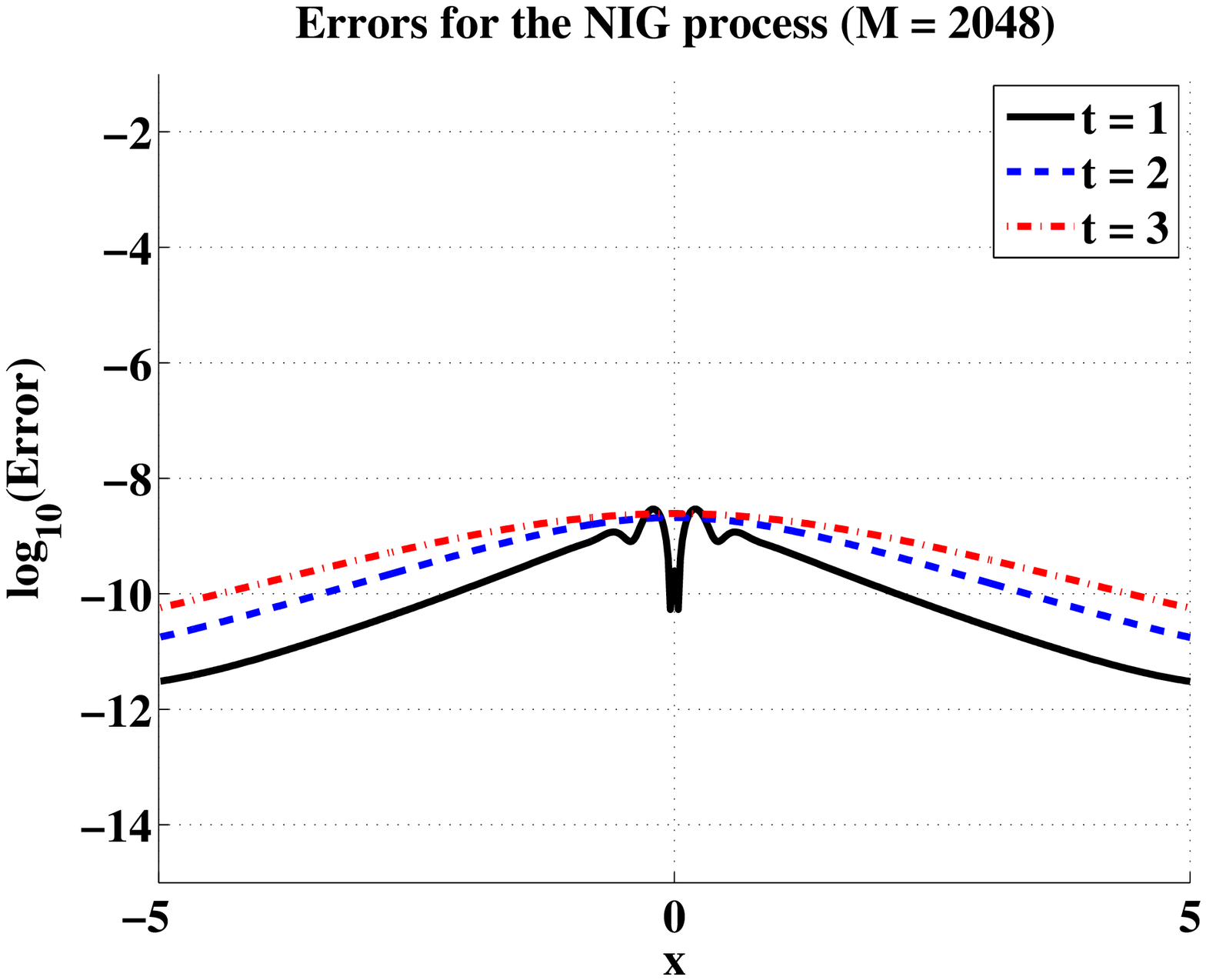}
\caption{Errors of the numerical solutions for Example~\ref{ex:NIG} 
when $t = 1,2,3$ and $M=2^{11}$ (i.e., $N=2^{8}$)}
\label{fig:err_NIG_M_2_11}
\end{minipage}

\end{center}
\end{figure}

\begin{figure}[H]
\begin{center}

\begin{minipage}[t]{0.45\linewidth}
\includegraphics[width = \linewidth]{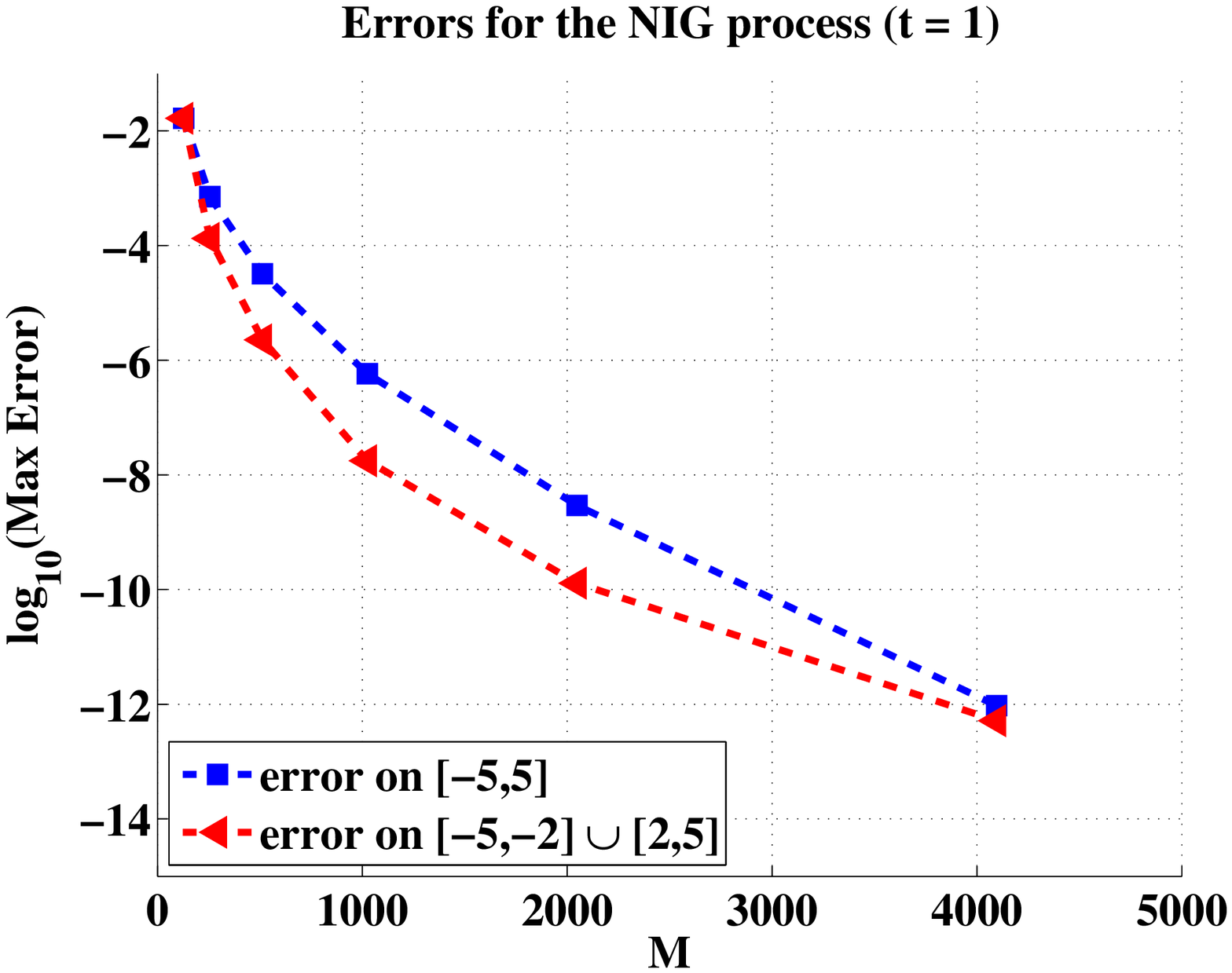}
\caption{Errors of the numerical solutions 
on $[-5, 5]$ and $[-5, -2] \cup [2, 5]$ for Example~\ref{ex:NIG} 
when $t = 1$ and $M = 2^{7}, \ldots,  2^{12}$ (i.e., $N = 2^{4}, \ldots,  2^{9}$)}
\label{fig:err_NIG_t1}
\end{minipage}\quad 
\begin{minipage}[t]{0.45\linewidth}
\includegraphics[width = \linewidth]{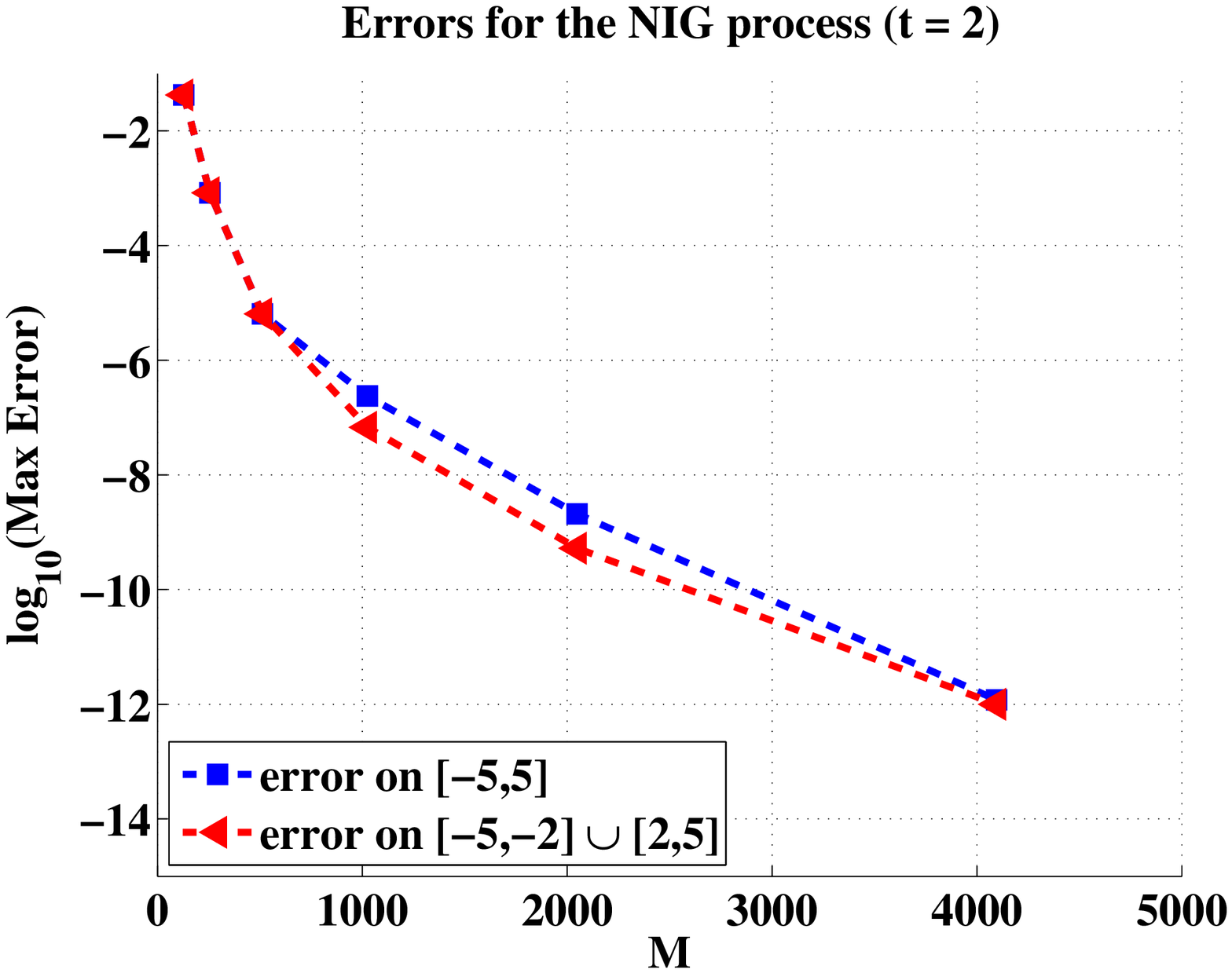}
\caption{Errors of the numerical solutions 
on $[-5, 5]$ and $[-5, -2] \cup [2, 5]$ for Example~\ref{ex:NIG} 
when $t = 2$ and $M = 2^{7}, \ldots,  2^{12}$ (i.e., $N = 2^{4}, \ldots,  2^{9}$)}
\label{fig:err_NIG_t2}
\end{minipage}

\end{center}
\end{figure}

\begin{figure}[H]
\begin{center}

\begin{minipage}[t]{0.45\linewidth}
\includegraphics[width = \linewidth]{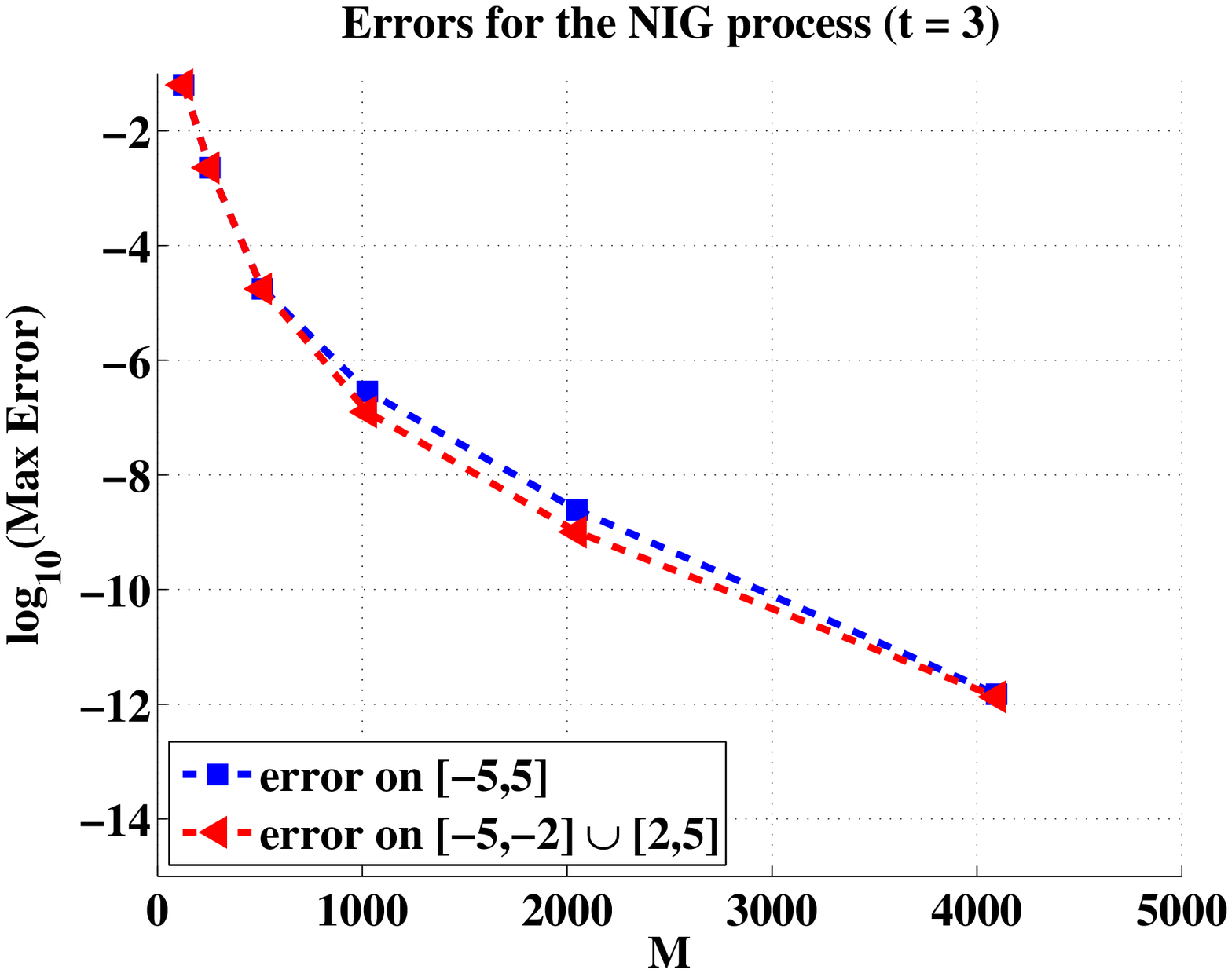}
\caption{Errors of the numerical solutions 
on $[-5, 5]$ and $[-5, -2] \cup [2, 5]$ for Example~\ref{ex:NIG} 
when $t = 3$ and $M = 2^{7}, \ldots,  2^{12}$ (i.e., $N = 2^{4}, \ldots,  2^{9}$)}
\label{fig:err_NIG_t3}
\end{minipage}\quad 
\begin{minipage}[t]{0.45\linewidth}
\includegraphics[width = \linewidth]{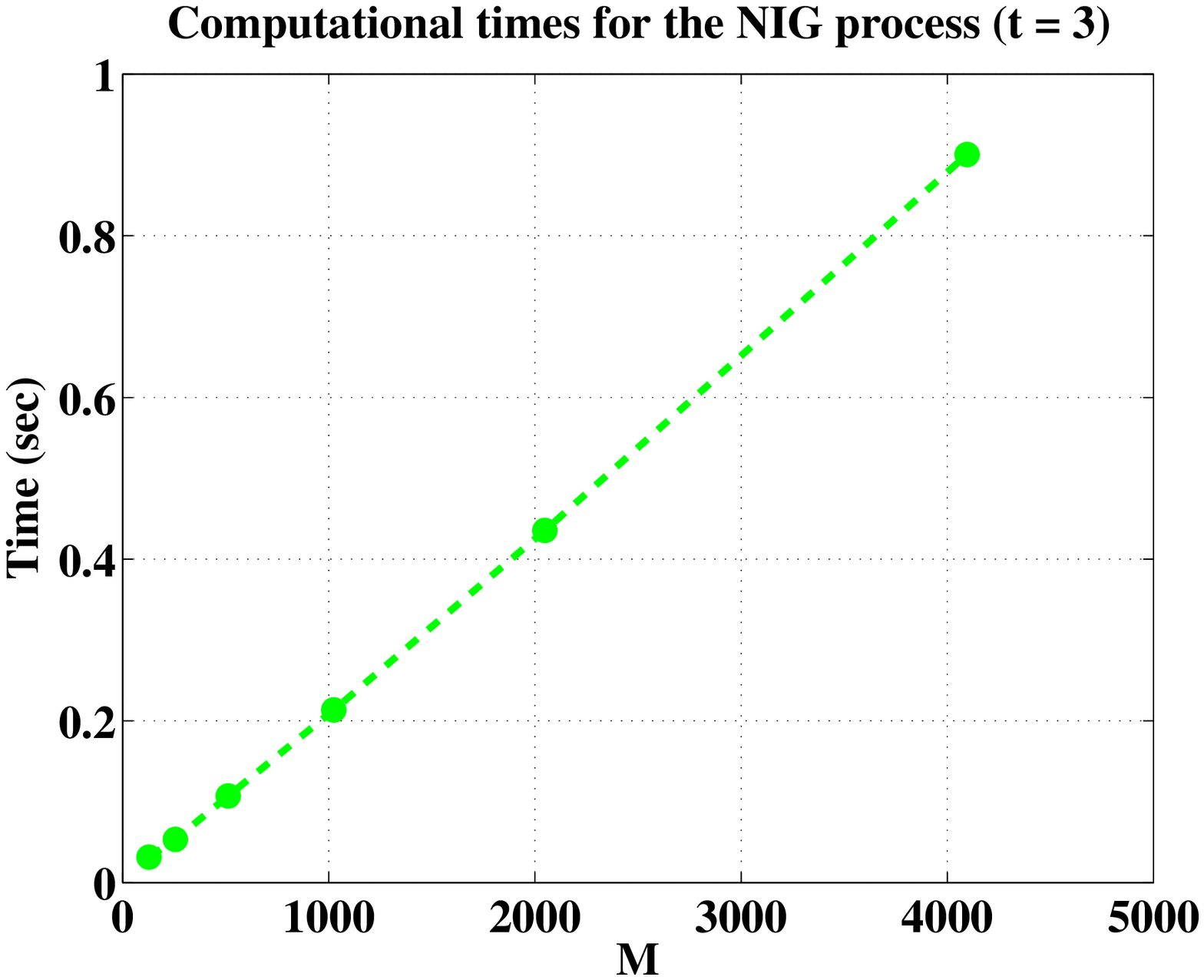}
\caption{Computational times for Example~\ref{ex:VG} when $t = 3$
and $M = 2^{7}, \ldots,  2^{12}$ (i.e., $N = 2^{4}, \ldots,  2^{9}$)}
\label{fig:time_NIG_t3}
\end{minipage}

\end{center}
\end{figure}

Results of these numerical experiments are shown below.
First, 
for reference, 
the exact solutions 
of~\eqref{eq:exact_sol_VG} of Example~\ref{ex:VG} 
for $t = 1,2,3$ 
are displayed in  
Figure~\ref{fig:sol_VG}. 
The errors of Example~\ref{ex:VG} for $M = 2^{11}$ 
and $t = 1,2,3$ are plotted in Figure~\ref{fig:err_VG_M_2_11}. 
The maximums of the errors of Example~\ref{ex:VG} 
on the intervals $[-x_{u}, x_{u}]$ and $[-x_{u}, -x_{l}] \cup [x_{l}, x_{u}]$
for every $M$'s are plotted for $t=1, 2$, and $3$
in Figures~\ref{fig:err_VG_t1}, \ref{fig:err_VG_t2}, and~\ref{fig:err_VG_t3}, 
respectively. 
The errors on $[-x_{u}, x_{u}]$ are computed 
to observe the errors which are not estimated by Theorem~\ref{thm:EulerFTError}.
The computational times of Example~\ref{ex:VG} for only $t=3$ are shown by Figure~\ref{fig:time_VG_t3}
because ones for $t = 1,2$ are considerably similar. 
Next, 
the exact solutions 
of~\eqref{eq:exact_sol_NIG} of Example~\ref{ex:NIG}
for $t = 1,2,3$ 
are displayed in Figure~\ref{fig:sol_NIG}.
The errors of Example~\ref{ex:NIG} for $M = 2^{11}$ 
and $t = 1,2,3$ are plotted in Figure~\ref{fig:err_NIG_M_2_11}. 
The maximums of the errors of Example~\ref{ex:NIG}
on the intervals $[-x_{u}, x_{u}]$ and $[-x_{u}, -x_{l}] \cup [x_{l}, x_{u}]$
for every $M$'s are plotted for $t=1, 2$, and $3$
in Figures~\ref{fig:err_NIG_t1}, \ref{fig:err_NIG_t2}, and~\ref{fig:err_NIG_t3}, 
respectively. 
The computational times of Example~\ref{ex:NIG} for only $t=3$ are shown by Figure~\ref{fig:time_NIG_t3}
because ones for $t = 1,2$ are considerably similar. 

The errors of Example~\ref{ex:VG} on the interval $[-x_{u}, -x_{l}] \cup [x_{l}, x_{u}]$
seems to have order $\mathrm{O}(\exp(-c \sqrt{M}))$ for some  $c>0$
according to Figures~\ref{fig:err_VG_t1}--\ref{fig:err_VG_t3}.  
This observation, Theorems~\ref{thm:SG} and~\ref{thm:EulerFTError}, and Remark~\ref{rem:ErrSGindef}
imply that the leading error occurs in Step~2 or~3 of the proposed method. 
On the other hand, 
in particular for $t=1$, 
the errors of Example~\ref{ex:VG} on the interval $[-x_{u}, x_{u}]$
are worse than ones on $[-x_{u}, -x_{l}] \cup [x_{l}, x_{u}]$. 
We can guess that this phenomenon is due to the cusp of the solution~\eqref{eq:exact_sol_VG}
at the origin shown by Figure~\ref{fig:sol_VG}. 
In fact, as time $t$ increases, 
the peakedness of the solution becomes gentler
and the errors around the origin improve. 
In addition, the computational times shown by Figure~\ref{fig:time_NIG_t3} 
are approximately consistent with the theoretical estimate 
$\mathrm{O}(N \log N) = \mathrm{O}(M \log M)$. 
As for the results of Example~\ref{ex:NIG}, 
we can obtain similar observations for the errors on the interval $[-x_{u}, -x_{l}] \cup [x_{l}, x_{u}]$
and the computational times. 
However, 
the errors of Example~\ref{ex:NIG} on the interval $[-x_{u}, x_{u}]$ 
are as good as the ones on $[-x_{u}, -x_{l}] \cup [x_{l}, x_{u}]$,
which may be because solution~\eqref{eq:exact_sol_NIG} 
does not have sharp cusp for $t = 1,2,3$. 

\section{Concluding remarks}
\label{sec:concl}

In this paper, 
we proposed a fast and accurate numerical method 
to solve the Kolmogorov forward equations~\eqref{eq:PIDE_LM}
of the scalar L\'evy processes with symmetric measures~\eqref{eq:Lv_measure}. 
The method consists of the three steps presented 
in Sections~\ref{sec:Outline} and~\ref{sec:NumScheme}. 
Step~1 and~3 are respectively based on 
accurate numerical formulas~\eqref{eq:DE_FT} and~\eqref{eq:LastEulerFFT}
for the Fourier transform proposed by \cite{bib:OouraEuler2001, bib:OouraDE-FT2005}, 
which are respectively combined with 
the nonuniform FFT and the fractional FFT to speed up the computations.
Step~2 requires numerical indefinite integration on the equispaced grids. 
This computation is performed using formula~\eqref{eq:SG_indefint_exchanged_plus} 
obtained by integrating the sinc-Gauss sampling formula~\eqref{eq:SG}
and combining the resultant convolution in~\eqref{eq:SG_indefint_exchanged_plus} with the FFT. 
The numerical solutions by the proposed method 
seemed to be exponentially convergent 
on the interval without sharp cusps of the corresponding exact solutions. 
Furthermore, the real computational times were approximately consistent 
with the theoretical estimate $\mathrm{O}(N \log N)$, 
where $N$ is the half of the number of the points $x$
on which the approximations of the solutions $p(x,t)$ were computed for a fixed $t$. 
As subjects of future works, 
we can consider the followings: 
the rigorous theoretical estimate of the errors of the proposed method, 
the optimal determination of the parameters based on the estimate, 
the comparison of the proposed method with other similar methods,  
and the extension of the method to broader class of L\'evy processes.

\section*{Acknowledgments}

The author would like to thank 
Prof.~L.~N.~Trefethen for his valuable comments regarding the 
sinc-Gauss indefinite integration formula in Step~2 of the proposed method 
in a private seminar at the University of Tokyo in March 2014. 
He also informed the author about reference~\cite{bib:HaleTownsend}. 
This work is supported by JSPS KAKENHI Grant Number 24760064.



\appendix

\section*{Appendix A.\ Computation of the integrals of the sinc-Gauss kernel}
\label{sec:SGint}

In this section, 
we propose an efficient method to compute the values of $G_{r}(k)$ in~\eqref{eq:def_SGIndefInt}, 
the integrals of the sinc-Gauss kernel. 
Let $F_{\mathrm{SG}}(\omega)$ be the Fourier transform of the sinc-Gauss kernel
\begin{align}
F_{\mathrm{SG}}(\omega) = \int_{-\infty}^{\infty} 
\left[ \mathop{\mathrm{sinc}}(x)\, \exp\left( - \frac{x^{2}}{2r^{2}} \right) \right] 
\exp(-\i\, \omega \, x)\, \d x.
\end{align}
Then, the function $F_{\mathrm{SG}}(\omega)$ is written in the form
\begin{align}
F_{\mathrm{SG}}(\omega) 
= 
\frac{1}{2} 
\left[
\mathop{\mathrm{erf}} \left( \frac{r (\omega + \pi)}{\sqrt{2}} \right)
-
\mathop{\mathrm{erf}} \left( \frac{r (\omega - \pi)}{\sqrt{2}} \right)
\right], 
\label{eq:SG_FT}
\end{align}
where $\mathop{\mathrm{erf}}$ is the error function defined as
\begin{align}
\mathop{\mathrm{erf}}(\xi) 
= \frac{2}{\sqrt{\pi}} \int_{-\infty}^{\xi} \exp(-t^{2})\, \mathrm{d}t. 
\end{align}
Using the function $F_{\mathrm{SG}}(\omega)$ in~\eqref{eq:SG_FT}, 
we have
\begin{align}
G_{r}(k+1) - G_{r}(k)
& =
\int_{k}^{k+1} 
\left(
\frac{1}{2 \pi} \int_{-\infty}^{\infty} F_{\mathrm{SG}}(\omega) \, 
\exp(\i\, x\, \omega )\, \d \omega 
\right)\, \d x \notag \\
& = 
\frac{1}{2 \pi} 
\int_{-\infty}^{\infty} 
F_{\mathrm{SG}}(\omega) 
\left(
\int_{k}^{k+1} 
\exp(\i\, x\, \omega )\, \d x 
\right) \, \d \omega \notag \\
& =
\frac{1}{2 \pi} 
\int_{-\infty}^{\infty} 
F_{\mathrm{SG}}(\omega) \, 
\mathop{\mathrm{sinc}}(\omega / (2 \pi))\, 
\exp(\i\, \omega / 2)\, 
\exp(\i\, k\, \omega )\, \d \omega, 
\label{eq:G_r_diff_IFT}
\end{align}
which is the inverse Fourier transform of the function 
$F_{\mathrm{SG}}(\omega) 
\mathop{\mathrm{sinc}}(\omega / (2 \pi)) 
\exp(\i\, \omega / 2)$. 
Since the function $F_{\mathrm{SG}}(\omega)$ rapidly decays as $|\omega| \to \infty$ on $\mathbf{R}$, 
applying the mid-point rule to integral~\eqref{eq:G_r_diff_IFT}, 
we can accurately approximate its values as 
\begin{align}
G_{r}(k+1) - G_{r}(k)
& \approx
\frac{h'}{2 \pi} 
\sum_{l' = -M+1}^{M} 
F_{\mathrm{SG}}(l' h')\, 
\mathop{\mathrm{sinc}}(l' h' / (2 \pi))\, 
\exp(\i\, l' h' / 2)\, 
\exp(\i\, k\, l' h'), 
\label{eq:G_r_diff_IDFT}
\end{align}
where $h' = 2\pi/M$.
This approximation is based on a similar principle as that of formula~\eqref{eq:LastEulerFFT}. 
Then, applying the fractional FFT to~\eqref{eq:G_r_diff_IDFT}, 
we can obtain the approximate values of $G_{r}(k+1) - G_{r}(k)$ for 
$k = 0, 1, \ldots, \lfloor M/2 \rfloor$ in $\mathrm{O}(M \log M)$ time. 
Finally, 
adding them sequentially from $k=0$ to $k = \lfloor M/2 \rfloor$,
we can compute the approximations of $G_{r}(k)$ for $k = 0, 1, \ldots, \lfloor M/2 \rfloor + 1$
in $\mathrm{O}(M)$ time. 


\begin{thebibliography}{10}

\bibitem[Applebaum(2009)Applebaum]
{bib:Apple_LevyText_2009}
\textsc{Applebaum, D.} (2009)
\newblock \emph{L\'evy Processes and Stochastic Calculus}, 2nd edn.
\newblock Cambridge: Cambridge University Press.

\bibitem[Bailey \& Swarztrauber(1991)Bailey \& Swarztrauber]
{bib:BailSwarFRFT1991}
\textsc{Bailey, D. H. \& Swarztrauber, P. N.} (1991) 
\newblock The fractional Fourier transform and applications.
\newblock \emph{SIAM Rev.}, \textbf{33}, 389--404.

\bibitem[Bueno-Orovio \emph{et al.}(2014)Bueno-Orovio, Kay \& Burrage]
{bib:Bueno_FourierFrac_2014}
\textsc{Bueno-Orovio, A., Kay, D. \& Burrage, K.} (2014)
\newblock Fourier spectral methods for fractional-in-space reaction-diffusion equations.
\newblock \emph{BIT Numer. Math.}, DOI 10.1007/s10543-014-0484-2.

\bibitem[Carr \& Madan(1999)Carr \& Madan]
{bib:CarrMadanPriceFFT1999}
\textsc{Carr, P. \& Madan, D.~B.} (1999)
\newblock Option valuation using the fast Fourier transform.  
\newblock \emph{J. Comput. Finance}, \textbf{2}, 61--73.

\bibitem[Chourdakis(2005)Chourdakis]
{bib:ChourdakisPriceFracFFT2005}
\textsc{Chourdakis, K.}  (2005)
\newblock Option pricing using the fractional FFT. 
\newblock \emph{J. Comput. Finance}, \textbf{8}, 1--18.

\bibitem[Cont \& Voltchkova(2005)Cont \& Voltchkova]
{bib:ContVoltchkova_PIDELevy_2005}
\textsc{Cont, R. \& Voltchkova, E.} (2005)
\newblock Integro-differential equations for option prices in exponential Levy models.
\newblock \emph{Finance Stochast.}, \textbf{9}, 299--325.

\bibitem[Duquesne \emph{et al.}(2010)Duquesne, Reichmann, Sato \& Schwab]
{bib:LevyMattersI_2010}
\textsc{Duquesne, T., Reichmann, O., Sato, K. \& Schwab, C.} (2010)
\newblock \emph{L\'evy Matters I:  
Recent Progress in Theory and Applications: Foundations, Trees and Numerical Issues in Finance}. 
Lecture Notes in Mathematics, vol. 2001.  
\newblock Heidelberg: Springer. 

\bibitem[Dutt \& Rokhlin(1993)Dutt \& Rokhlin]
{bib:DuttRokhlin_NFFT_1993}
\textsc{Dutt, A. \& Rokhlin, V.} (1993) 
\newblock Fast Fourier transforms for nonequispaced data. 
\newblock \emph{SIAM J. Sci. Comput.}, \textbf{14}, 1368--1393.

\bibitem[Dutt \& Rokhlin(1995)Dutt \& Rokhlin]
{bib:DuttRokhlin_NFFT_1995}
\textsc{Dutt, A. \& Rokhlin, V.} (1995)
\newblock Fast Fourier transforms for nonequispaced data II. 
\newblock \emph{Appl.~Comput.~Harmon.~Anal.}, \textbf{2}, 85--100.

\bibitem[Fang \& Oosterlee(2008)Fang \& Oosterlee]
{bib:FangOost_COS_2008}
\textsc{Fang, F \& Oosterlee, C. W.} (2008)
\newblock A novel pricing method for European options based on Fourier-cosine series expansions.
\newblock \emph{SIAM J.~Sci.~Comput.}, \textbf{31}, 826--848.

\bibitem[Gao \emph{et al.}(2013)Gao, Duan \& Li]
{bib:Gao_FPE_symLevy_2013}
\textsc{Gao, T., Duan, J. \& Li, X.} (2013)
\newblock Fokker-Planck equations for stochastic dynamical systems with symmetric L\'evy motions. 
\newblock  arXiv:1310.7677.

\bibitem[Gardiner(2009)Gardiner]
{bib:Gardiner_StocMeth_2009}
\textsc{Gardiner, C.} (2009)
\newblock \emph{Stochastic Methods: A Handbook for the Natural and Social Sciences}, 4th edn. 
\newblock Heidelberg: Springer. 

\bibitem[Garreau \& Kopriva(2013)Garreau \& Kopriva]
{bib:GarreauKopriva2013}
\textsc{Garreau, P. \& Kopriva, D.} (2013)
\newblock A spectral element framework for option pricing under general exponential L\'evy processes.
\newblock \emph{J. Sci. Comput.}, \textbf{57}, 390--413. 

\bibitem[Greengard \& Lee(2004)Greengard \& Lee]
{bib:GreengardLee_NFFT_2004}
\textsc{Greengard, L. \& Lee, J. Y.} (2004)
\newblock Accelerating the nonuniform fast Fourier transform.
\newblock \emph{SIAM Rev.}, \textbf{46}, 443--454.

\bibitem[Hale \& Townsend(2014)Hale \& Townsend]
{bib:HaleTownsend}
\textsc{Hale, N. \& Townsend, A.} (2014)
\newblock An algorithm for the convolution of Legendre series. 
\newblock \emph{SIAM J. Sci. Comput.}, \textbf{36}, A1207--A1220. 

\bibitem[Huang \emph{et al.}(2014)Huang, Nie \& Tang]
{bib:Huang_FDSpec_2014}
\textsc{Huang, J., Nie, N. \& Tang, Y.} (2014)
\newblock A second order finite difference-spectral method for space fractional diffusion equations. 
\newblock \emph{Sci. China Math.}, \textbf{57}, 1303--1317.

\bibitem[Huang \& Oberman(2013)Huang \& Oberman]
{bib:Huang_FracLaplaceI_2013}
\textsc{Huang, Y. \& Oberman, A.} (2013)
\newblock Numerical methods for the fractional Laplacian Part I: a finite difference-quadrature approach. 
\newblock arXiv:1311.7691. 

\bibitem[Kozubowski \emph{et al.}(2006)Kozubowski, Meerschaert \& Podg\'orski]
{bib:Kozu_etal_FLM_2006}
\textsc{Kozubowski, T. J., Meerschaert, M. M. \& Podg\'orski, K.} (2006)
\newblock Fractional Laplace motion,  
\newblock \emph{Adv. Appl. Prob.}, \textbf{38}, 451--464. 

\bibitem[Kwok \emph{et al.}(2012)Kwok, Leung \& Wong]
{bib:Kwok_etal_2012}
\textsc{Kwok, Y. K., Leung, K. S. \& Wong, H. Y.} (2012)
\newblock Efficient options pricing using the fast Fourier transform. 
\emph{Handbook of Computational Finance}
(J.-C. Duan \emph{et al}. eds). Berlin: Springer, pp. 579--604.

\bibitem[Lee \emph{et al.}(2012)Lee, Liu \& Sun]
{bib:Lee_etal_FastExpIntOP_2012}
\textsc{Lee, S. T., Liu, X. \& Sun, H.-W.} (2012)
\newblock Fast exponential time integration scheme for option pricing with jumps. 
\newblock \emph{Numer. Linear Algebra Appl.}, \textbf{19}, 87--101.

\bibitem[Li \emph{et al.}(2012)Li, Deng \& Wu]
{bib:Li_FDM_FPE_2012}
\textsc{Li, C., Deng, W. \& Wu, Y.} (2012)
\newblock Finite difference approximations and dynamics simulations for the Levy fractional Klein-Kramers equation. 
\newblock \emph{Numer. Methods Partial Differential Equations}, \textbf{28}, 1944--1965. 

\bibitem[Lenzi \emph{et al.}(2003)Lenzi, Mendes, Kwok \& Malacarne]
{bib:Lenzi_etal_FracFPE_2003}
\textsc{Lenzi, E. K., Mendes, R. S., Kwok, S. F. \& Malacarne, L. C.} (2003)
\newblock Anomalous diffusion: fractional Fokker-Planck equation and its solutions. 
\newblock \emph{J.~Math.~Phys.}, \textbf{44}, 2179--2185.

\bibitem[Meerschaert(2004)Meerschaert \& Tadjeran]
{bib:Meerschaert_FDM_Flow_2004}
\textsc{Meerschaert. M. M. \& Tadjeran, C.} (2004)
\newblock Finite difference approximations for fractional advection-dispersion flow equations.
\newblock \emph{J. Comput. Appl. Math.}, \textbf{172}, 65--77.

\bibitem[Ooura(2001)Ooura]
{bib:OouraEuler2001}
\textsc{Ooura, T.} (2001) 
\newblock A continuous Euler transformation and its application to 
Fourier transform of a slowly decaying function. 
\newblock \emph{J. Compt. Appl. Math.}, \textbf{130}, 259--270.

\bibitem[Ooura(2005)Ooura]
{bib:OouraDE-FT2005}
\textsc{Ooura, T.} (2005)
\newblock A double exponential formula for the Fourier transforms. 
\newblock \emph{Publ. RIMS Kyoto Univ.}, \textbf{41}, 971--977.

\bibitem[Potts \emph{et al.}(2001)Potts, Steidl \& Tasche]
{bib:Potts_etal_NFFT_2001}
\textsc{Potts, D., Steidl, G. \& Tasche, M.} (2001)
\newblock Fast Fourier transforms for nonequispaced data: A tutorial. 
\emph{Modern Sampling Theory: Mathematics and Applications}.
(J. J. Benedetto \& P. Ferreira, eds).
\newblock Boston: Birkh\"auser, pp. 249--274.

\bibitem[Sabatier \emph{et al.}(2007)Sabatier, Agrawal \& Tenreiro Machado]
{bib:Sabatier_etal_FracCal_2007}
\textsc{Sabatier, J., Agrawal, O. P. \& Tenreiro Machado, J. A.} (2007)
\newblock \emph{Advances in Fractional Calculus: Theoretical Developments and Applications in Physics and Engineering}. 
\newblock Springer.

\bibitem[Steidl(1998)Steidl]
{bib:Steidl_NFFT_1998}
\textsc{Steidl, G.} (1998)
\newblock A note on fast Fourier transforms for nonequispaced grids.
\newblock \emph{Adv. Comput. Math.}, \textbf{9}, 337--352.

\bibitem[Tanaka(2014a)Tanaka]
{bib:KTanaka_EulerFFT_2014}
\textsc{Tanaka, K.} (2014a) 
\newblock Error control of a numerical formula for the Fourier transform by Ooura's continuous Euler transform and fractional FFT. 
\newblock \emph{J. Comput. Appl. Math.}, \textbf{266}, 73--86.

\bibitem[Tanaka(2014b)Tanaka]
{bib:KTanaka_Matlab_2014}
\textsc{Tanaka, K.} (2014b) 
\newblock Matlab codes for the symmetric Levy processes. \\
\newblock \url{https://github.com/KeTanakaN/mat_symLevy_FT_codes}\  
(accessed 1 August 2014). 

\bibitem[Tanaka \emph{et al.}(2008)Tanaka, Sugihara \& Murota]
{bib:KTanaka_etal_SG_2008}
\textsc{Tanaka, K., Sugihara, M. \& Murota, K.} (2008)
\newblock Complex-analytic approach to the sinc-Gauss sampling formula. 
\newblock \emph{Japan J. Indust. Appl. Math.}, \textbf{25}, 209--231.

\bibitem[Tanaka \emph{et al.}(2009)Tanaka, Sugihara, Murota, \& Mori]
{bib:KTanaka_etal_DE_2009}
\textsc{Tanaka, K., Sugihara, M., Murota, K. \& Mori, M.} (2009)
\newblock Function classes for double exponential integration formulas. 
\newblock \emph{Numer. Math.}, \textbf{111}, 631--655.

\bibitem[Yan(2013)Yan]
{bib:Yan_FFPE_Laplace_2013}
\textsc{Yan, L.} (2013) 
\newblock Numerical solutions of fractional Fokker-Planck equations
using iterative Laplace transform method. 
\emph{Abstr. Appl. Anal.}, \textbf{2013}, Art. ID 465160.

\bibitem[Zhao \& Lib(2012)Zhao \& Lib]
{bib:Zhao_NumFPE_2012}
\textsc{Zhao, Z. \& Lib, C.} (2012)
\newblock A numerical approach to the generalized nonlinear fractional Fokker-Planck equation. 
\newblock \emph{Comput. Math. Appl.}, \textbf{64}, 3075--3089. 

\end{thebibliography}
\end{document}